\documentclass{article}
\usepackage{amssymb,amsmath,amsthm,verbatim,IEEEtrantools,macros,mathtools}
\usepackage[margin=1in]{geometry}
\newcommand{\Exact}{\operatorname{Exact}}
\newcommand{\supp}{\operatorname{supp}}
\newcommand{\cyl}{\operatorname{cyl}}
\begin{document}

\title{Fourier Dimension Estimates for Sets of Exact Approximation Order: The Badly-Approximable Case}

\author{Robert Fraser and Reuben Wheeler}

\maketitle

\begin{abstract}
We show for decreasing, positive approximation functions $\psi$ such that $\tau = \lim_{q \to \infty} \frac{\log \psi(q)}{\log q} < \frac{13 + \sqrt{73}}{8}$ and such that $q^2 \psi(q) \to 0$ that the set $\Exact(\psi)$ of numbers approximable to the exact order $\psi$ has positive Fourier dimension. This implies that the set $\Exact(\psi)$ contains normal numbers.
\end{abstract}
\section{Introduction}
\subsection{Metrical Diophantine Approximation}
The study of Diophantine approximation concerns approximation of real numbers by rational numbers with small denominators. The classical Dirichlet principle states that if $x \in \mathbb{R}$, then there exist infinitely many rational numbers $\frac{p}{q}$ such that
\[\left|x - \frac{p}{q} \right| \leq q^{-2}.\]
For some irrational numbers, the exponent in the upper bound $q^{-2}$ can be significantly improved; for others, $q^{-2}$ is optimal up to a multiplicative constant. These sets of numbers are called the \textbf{well-approximable numbers} and \textbf{badly-approximable numbers}, respectively.

Let $\tau > 2$ be a real number. We define the set $E(\tau)$ of \textbf{well-approximable numbers} to consist of those $x \in \mathbb{R}$ satisfying

\[\left|x - \frac{p}{q} \right| \leq q^{-\tau} \quad \text{for infinitely many pairs of integers $(p,q)$.}\]

On the other hand, let $N \geq 2$ be a positive integer. We define the set $\textbf{Bad}(N)$ of \textbf{badly-approximable numbers} to be the set of those real numbers $x$ such that every partial quotient in the continued fraction expansion of $x$ is bounded above by $N$. If $x \in \textbf{Bad}(N)$, then $x$ satisfies the estimate
\[\left|x - \frac{p}{q} \right| \geq \frac{1}{2(N + 1)q^2}\]
for every rational number $\frac{p}{q}$; conversely, if $x$ satisfies the inequality
\[\left|x - \frac{p}{q} \right| \geq c q^{-2}\] 
for every rational number $\frac{p}{q}$, then $x \in \textbf{Bad}(N)$ where $N = \lceil 2 c^{-2} \rceil$. 

More generally, given a finite set $S \subset \mathbb{N}$ with at least two elements, we define the set $\textbf{Bad}'(S)$ to be the set of numbers with continued fractions whose partial quotients lie in $S$.

\textit{Metrical Diophantine approximation} is the study of metrical properties of sets such as $E(\tau)$ and $\textbf{Bad}(N)$ arising from Diophantine approximation. One metrical property that is of great interest is Hausdorff dimension. The Hausdorff dimension of $E(\tau)$ is computed by Jarník \cite{Jarnik29} and Besicovitch \cite{Besicovitch34} to be $\frac{2}{\tau}$. Good \cite{Good41} obtains an implicit description of the Hausdorff dimension of $\textbf{Bad}(N)$.

It is also interesting to consider the Hausdorff dimension of sets of numbers satisfying more precise approximation conditions than those of $E(\tau)$ or $\text{Bad}(N)$. For example, suppose $q \psi(q)$ is a decreasing function. For example, a result of Khintchine \cite{Khintchine24} states that the set $E(\psi)$ of numbers $x$ satisfying the bound
\[\left|x - \frac{p}{q} \right| \leq \psi(q) \quad \text{for infinitely many pairs of integers $(p,q)$}\]
will be Lebesgue-null if and only if the sum of $q \psi(q)$ is finite. 

Another direction in the study of metrical Diophantine approximation involves considering analogues of the well-approximable numbers in higher dimensions. In particular, Bovey and Dodson \cite{BoveyDodson86} compute the Hausdorff dimension of higher-dimensional analogues of the well-approximable numbers.  See the book of Bernik and Dodson \cite{BernikDodson99} for a survey of metrical results in Diophantine approximation.

More generally, the well-approximable numbers are an example of a \textbf{lim sup set}. A lim sup set is a set that can be expressed in the form
\[\bigcap_{n=1}^{\infty} \bigcup_{k = n}^{\infty} A_n\]
for some balls $A_n$. A mass transference principle of Beresnevich and Velani \cite{BeresnevichVelani15} can be used to estimate the Hausdorff dimension of any lim sup set by estimating the Lebesgue measure of an appropriate ``dilated" version of the set. This mass transference principle can be used to recover all of the above results on the Hausdorff dimension of the well-approximable numbers, but cannot be used to obtain results on the badly-approximable numbers, which are not a lim sup set.
\subsection{Normal numbers}
Let $a \geq 2$ be an integer and let $x$ be a real number. We say that the number $x$ is \textbf{simply normal to base $a$} if each digit in the base-$a$ expansion of $x$ occurs with the same frequency. The number $x$ is said to be \textbf{normal to base $a$} if for each integer $l$, each $l$-tuple of digits occurs with the same frequency in the base-$a$ expansion of $x$. It is easy to see that $x$ is normal to base $a$ if and only if the fractional parts $\{a^n x\}_{n=1}^{\infty}$ are uniformly distributed modulo $1$. A number is said to be \textbf{normal} if it is normal to base $a$ for every $a \geq 2$. 

The problem of determining whether a specific irrational number is normal is in general a very difficult one. For example, the problem of determining whether numbers such as $e$ or $\pi$ are normal is still open. On the other hand, a simple argument shows that Lebesgue-almost-all numbers are normal. In particular, any subset of $\mathbb{R}$ of positive Lebesgue measure must contain normal numbers.

As a compromise between the trivial problem of locating normal numbers in subsets of $\mathbb{R}$ of positive Lebesgue measure and the intractable problem of determining whether a specific irrational number is normal, it is natural to ask for conditions that can guarantee that a Lebesgue-null subset $E$ of $\mathbb{R}$ contains normal numbers. It is not enough to assume that $E$ has positive Hausdorff dimension. For example, it is simple to see that the Cantor middle-thirds set does not even contain any numbers that are simply normal to base $3$.

A condition on $E$ that guarantees the existence of normal numbers is given by Davenport, Erdős, and Leveque. This condition states that if $\mu$ is a Borel probability measure, then the sequence $\{a^n x\}_{n=1}^{\infty}$ is uniformly distributed modulo $1$ for $\mu$-a.e. $x$ if 
\begin{equation}\label{DELcondition1}
\sum_{N=1}^{\infty} N^{-3} \sum_{j=1}^N \sum_{k=1}^N \hat \mu(m(a^j - a^k)) < \infty.
\end{equation}
A simple calculation (see Pollington, Velani, Zafeiropoulos, and Zorin \cite{PollingtonVelaniZafeiropoulosZorin22}) shows that the condition \eqref{DELcondition1} will hold for all $a$ provided that $\mu$ satisfies the quantitative Fourier decay condition
\begin{equation}\label{DELcondition2}
\hat \mu(\xi) \lesssim C (\log \log \xi)^{-(1 + \epsilon)}
\end{equation}
for some $\epsilon > 0$.
\subsection{Fourier dimension and Hausdorff dimension}
We begin by presenting some classical results about Hausdorff and Fourier dimension. We recommend Mattila \cite{Mattila15} as a reference for the classical material discussed in this section.

In geometric measure theory, Frostman's lemma \cite{Frostman35} states that for any compact set (or, more generally, for any Borel set) $E \subset \mathbb{R}^n$, the Hausdorff dimension of $E$ is the supremum of all of those $0 \leq s \leq n$ such that there exists a Borel probability measure $\mu_s$ supported on $E$ satisfying the ball condition
\begin{equation}\label{ballConditionGeneral}
\mu(B) \leq C_s \operatorname{radius}(B)^s \quad \text{for all balls $B \subset \mathbb{R}^n$.}
\end{equation}
Intuitively, this condition states that $E$ has large Hausdorff dimension if it supports a measure that is ``spread out." 

A related condition is the energy integral estimate 
\begin{equation}\label{energyIntegral}
\mathcal{E}(\mu) := \iint_{\substack{x \in \mathbb{R}^n \\ y \in \mathbb{R}^n}} |x - y|^{-s} \, d \mu(x) \, d \mu(y) < \infty.
\end{equation}
The condition \eqref{energyIntegral} is almost equivalent to the condition \eqref{ballConditionGeneral}. If either of these conditions holds for some $s_0$, the other will hold for all $s < s_0$. Hence the Hausdorff dimension of a compact set (or more generally, a Borel set) $E$ is the supremum of those values of $s$ for which $E$ supports a measure satisfying the condition \eqref{energyIntegral}.

A further characterization of the Hausdorff dimension is obtained by obtaining a Fourier-analytic expression for the energy integral. In fact, an argument based on the Fourier transform of the Riesz kernel shows that there exists a constant $C(s,n)$ such that 
\begin{equation}\label{fourierEnergy}
\mathcal{E}(\mu) = C(s,n) \int_{\xi \in \mathbb{R}^n} |\hat \mu(\xi)|^2 |\xi|^{s - n} \, d \xi
\end{equation}
Hence the Hausdorff dimension of $E$ is the supremum over those values of $s$ for which $E$ supports a measure such that the integral \eqref{fourierEnergy} is finite.

The condition \eqref{fourierEnergy} is an $L^2$-averaged condition on the Fourier decay of the measure $\mu$. Up to an $\epsilon$-loss in the exponent $s$, the integral in \eqref{fourierEnergy} is finite if and only if $|\hat \mu(\xi)|$ decays like $|\xi|^{-s/2}$ ``on average". However, the condition \eqref{fourierEnergy} does not imply \textbf{any} pointwise bound on $\hat \mu(\xi)$, not even a bound on the order of \eqref{DELcondition2}. It is easy to see in dimensions $n \geq 2$ that a measure supported a hyperplane $H$ will have no Fourier decay in directions orthogonal to $H$. When $n = 1$, it can be shown that no Borel probability measure supported on the Cantor middle thirds set decays pointwise as $\xi \to \infty$. In fact, for Cantor middle-$\alpha$ sets $C_{\alpha}$, the pointwise Fourier decay of measures supported on $C_{\alpha}$ is related to the algebraic properties of $\alpha$; an excellent reference for this topic is Bugeaud \cite{Bugeaud12}.

It is therefore useful to introduce a notion of dimension that captures the \textit{pointwise} Fourier decay of measures supported on $E$. For a compact set (or Borel set) $E$, the \textit{Fourier dimension} of $E$ is defined to be the supremum over those values $0 \leq s \leq n$ such that $E$ supports a Borel probability measure $\mu_s$ satisfying the pointwise Fourier decay condition
\begin{equation}\label{fourierDimension}
|\hat \mu_s(\xi)| \leq C_s (1 + |\xi|)^{-s/2}.
\end{equation}
Observe that the condition \eqref{fourierDimension} implies the condition \eqref{DELcondition2}. Hence \textit{any compact subset (or Borel subset) of $\mathbb{R}$ of positive Fourier dimension contains normal numbers}.

For hypersurfaces in $\mathbb{R}^n$, the Fourier dimension of a set captures the \textit{curvature} of the set. For example, a sphere or paraboloid in $\mathbb{R}^n$ will have Fourier dimension $n - 1$, while a subset of a hyperplane will always have Fourier dimension $0$.

For fractal subsets of $\mathbb{R}^n$, the Fourier dimension captures the \textit{additive irregularity} of the set. So ``regular" sets like the Cantor middle-thirds set will have Fourier dimension $0$, but random Cantor sets will almost surely have positive Fourier dimension.

Fourier dimension is important in harmonic analysis because sets with positive Fourier dimension satisfy a \textbf{restriction estimate}. If $\mu$ is a Borel probability measure supported on a compact set $E \subset \mathbb{R}^n$, define the \textit{extension operator} $\mathcal{R}_{\mu}^*$ on $L^{\infty}(E)$ by the integral
\begin{equation}\label{restrictionOperator}
\mathcal{R}_{\mu}^*f(x) := \int_E e^{2 \pi i x \cdot \xi} f(\xi) \, d \xi.
\end{equation}
A classical result of Stein \cite{Stein67} and Tomas \cite{Tomas75} states that if $E$ is a hypersurface of $\mathbb{R}^n$ of nonvanishing Gaussian curvature, then the surface measure on $E$ satisfies the restriction estimate
\begin{equation}\label{restrictionEstimate}
\norm{\mathcal{R}_{\mu}^*f}_{L^p} \leq \norm{f}_{L^2},
\end{equation}
where the exponent $p \geq \frac{2n + 2}{n - 1}$. The value $\frac{2n + 2}{n - 1}$ cannot be improved. It was later shown by Mockenhaupt \cite{Mockenhaupt00} and Mitsis \cite{Mitsis02} that if $\mu$ is a measure satisfying the ball condition \eqref{ballConditionGeneral} with exponent $\alpha$ and satisfying the pointwise Fourier decay condition \eqref{fourierDimension} with exponent $\beta$, then $\mu$ satisfies the restriction estimate \eqref{restrictionEstimate} for any $p > \frac{4n - 4 \alpha + 2 \beta}{\beta}$. The corresponding restriction estimate at the endpoint $p = \frac{4n - 4 \alpha + 2 \beta}{\beta}$ is shown by Bak and Seeger \cite{BakSeeger11}. Observe that in the case of the surface measure of a surface with nonvanishing Gaussian curvature, we can take $\alpha = \beta = n - 1$ and recover the Stein-Tomas exponent $\frac{2n + 2}{n - 1}$.

Two classical results in the theory of Fourier dimension are due to Kaufman \cite{Kaufman80}, \cite{Kaufman81}. Kaufman computes the Fourier dimension of the well-approximable numbers $E(\tau)$ to be exactly $\frac{2}{\tau}$ \cite{Kaufman81}; this result is exposited in detail in a work of Bluhm \cite{Bluhm98}, and some Fourier dimension estimates for $n \geq 1$ are given by Hambrook \cite{Hambrook19}. All of these arguments are fundamentally based on the cancellation of the exponential sum $\sum_{p=0}^{q-1} e^{2 \pi i \frac{p s}{q}}$ for any integer $s$.

On the other hand, Kaufman \cite{Kaufman80} also provides a lower Fourier dimension estimate for the Fourier dimension of the Badly-approximable numbers with an argument based on the continued fraction expansion. Specifically, Kaufman shows that if $S$ is such that $\textbf{Bad}'(S)$ has Hausdorff dimension greater than $2/3$, then $\textbf{Bad}'(S)$ has positive Fourier dimension. A summary of many of the ideas from this paper will appear in Sections \ref{sec:oscArg}--\ref{sec:combEst}. The value $2/3$ from Kaufman's result has been improved to $1/2$ by Queffélec and Ramaré. Later, Hochman and Shmerkin \cite{HochmanShmerkin15} show that for any set $S \subset \mathbb{N}$ with at least two elements, the set $\textbf{Bad}'(S)$ will always contain normal numbers, and Sahlsten and Stevens \cite{SahlstenStevens20} have shown that for such $S$, $\textbf{Bad}'(S)$ will always have positive Fourier dimension.
\subsection{Sets of Exact Approximation Order}
Bugeaud \cite{Bugeaud03} introduces a set that shares some of the properties of the well-approximable numbers and some of the properties of the badly approximable numbers. This set is called the \textbf{set of numbers of exact approximation order}. Let $\psi : \mathbb{N} \to \mathbb{R}_{> 0}$ be a function. Bugeaud defines the set $\Exact(\psi)$ to be the set of numbers satisfying the pair of conditions
\begin{IEEEeqnarray*}{rCll}
\left|x - \frac{p}{q} \right| & \leq & \psi(q) \quad & \text{for infinitely many pairs of integers $(p,q)$} \\
\left|x - \frac{p}{q} \right| & \geq & (1 - c) \psi(q) \quad & \text{for any $c > 0$ and all $q > Q(c)$}.
\end{IEEEeqnarray*}
The set $\Exact(\psi)$ consists of numbers that can be approximated to the order $\psi$, but such that the approximation function $\psi$ cannot be improved by even a multiplicative constant. Hence $\Exact(\psi)$ will always be contained in $E(\psi)$. Although it may seem that $\Exact(\psi)$ should be much smaller than $E(\psi)$, Bugeaud \cite{Bugeaud03} shows that as long as $q^2 \psi(q)$ is decreasing, these two sets have the same Hausdorff dimension. In fact, Bugeaud estimates the Hausdorff dimension of $\Exact(\psi)$ by showing that the set $\Exact(\psi)$ contains a family of numbers that have continued fraction expansions with the property that most of the partial quotients grow very slowly, except for a thin subsequence of the partial quotients that grow at a rate close to $\frac{1}{q^2 \psi(q)}$. Hence it seems likely that a version of Kaufman's argument \cite{Kaufman80} for the badly-approximable numbers should be applicable to this set.

In a previous work \cite{FraserWheeler22}, we establish the following result by modifying Kaufman's construction for the \textit{well-approximable} numbers.
\begin{mythm}\label{largetau}
Let $\psi : \mathbb{N} \to (0, \infty)$ be a positive, decreasing function such that the limit
\[\tau := - \lim_{q \to \infty} \frac{\log \psi(q)}{\log q}\]
exists and is finite. Suppose further that $\tau > \frac{3 + \sqrt{5}}{2}$. Then $\Exact(\psi)$ has positive Fourier dimension.
\end{mythm}
In the work \cite{FraserWheeler22}, we explicitly obtain a lower bound on the Fourier dimension of $\Exact(\psi)$. Moreover, we establish a version of Theorem \ref{largetau} for inhomogeneous approximation. However, Theorem \ref{largetau} does not say anything for smaller values of $\tau$.

In this work we will obtain a Fourier dimension estimate for $\Exact(\psi)$ given suitable conditions on the approximation function $\psi$.

\begin{mythm}\label{mainthm}
Let $\psi: \mathbb{N} \to \mathbb{R}^+$ be any function satisfying the conditions
\begin{equation}\label{q2psiqbound}
0 \leq  q^2 \psi(q) \leq 1 \quad \text{for all $q$} \\
\end{equation}
\begin{equation}\label{q2psiqlimit}
q^2 \psi(q) \to 0 \quad \text{as  $q \to \infty$}
\end{equation}
\begin{equation}\label{psiorder}
\lim_{q \to \infty} - \frac{\log \psi(q)}{\log q} := \tau < \frac{13 + \sqrt{73}}{8}.
\end{equation}

Then $\Exact(\psi)$ has positive Fourier dimension.
\end{mythm}
\begin{myrmk}\label{mainthmrmk}
It is possible that the methods in this paper can be improved to obtain a larger range of $\tau$, as suggested by the extension of Kaufman's work by Queffélec and Ramaré \cite{QueffelecRamare03}. Nevertheless, the result is satisfactory because $\frac{13 + \sqrt{73}}{8}$ is larger than the value $\frac{3 + \sqrt{5}}{2}$ appearing in Theorem \ref{largetau}.
\end{myrmk}

\begin{myrmk}\label{mainthmrmk2}
The techniques in the proof of Theorem \ref{mainthm} do not seem to easily generalize to the inhomogeneous setting. In fact, the problem of showing that the inhomogeneous badly approximable numbers have positive Fourier dimension has been proposed by Pollington, Velani, Zafeiropoulos, and Zorin \cite{PollingtonVelaniZafeiropoulosZorin22}.
\end{myrmk}

\begin{myrmk}We obtain a small lower bound on Fourier dimension for a given exact order set. Improved error estimates for the measure approximation of Section \ref{sec:measAppr} or new arguments would be necessary to establish stronger Fourier dimension estimates.\end{myrmk}

\section{Notation}
\setcounter{subsection}{1}
For a finite sequence $G$ and a sequence $F$ (possibly infinite), we write $G \cdot F$ for their concatenation. For a finite sequence of natural numbers $G$, we write $G^*$ for the set $\lbrace G\rbrace \times \mathbb{N}^{\mathbb{N}}$, the set of sequences of natural numbers beginning with $G$. 

For real numbers $x$, we write $e(x)=e^{2\pi i x}$.

Given a finite sequence of natural numbers $G = (c_0, c_1, \ldots, c_t)$, we will use $K(G)$ to denote the denominator of the finite continued fraction $[c_0; c_1, \ldots, c_t]$, and we will use $K'(G)$ to denote the denominator, or continuant, of the finite continued fraction $[c_0; c_1, \ldots, c_{t-1}]$. Furthermore, we denote by $\cyl(G)$ the set of real numbers that have $[c_0;c_1,\ldots,c_t]$ as a convergent. More specifically, $\cyl(G)=[[c_0;c_1,\ldots,c_t],[c_0;c_1,\ldots,c_t+1])$ if $t$ is even and $\cyl(G)=([c_0;c_1,\ldots,c_t+1],[c_0;c_1,\ldots,c_t]]$ if $t$ is odd.

With $\psi$ as in Theorem \ref{mainthm}, it is convenient to define the function
\begin{equation}\label{rhodef}
\rho(q) := \frac{1}{q^2 \psi(q)}.
\end{equation}
Note that $\rho$ will satisfies the following:
\begin{equation}\label{rhoorder}
\lim_{q \to \infty} \frac{\log \rho(q)}{\log q} = \tau - 2 \leq \frac{-3 + \sqrt{73}}{8} 
\end{equation}
and
\begin{equation}\label{rholimit}
\lim_{q \to \infty} \rho(q) = \infty.
\end{equation}


Throughout this proof, we will have quantities $\epsilon$ and $N$ chosen so that $1 - \epsilon < \text{dim}_H \text{Bad}(N) < 1 - \frac{\epsilon}{10}$. Whenever a statement is made for ``sufficiently large $N$," it should be assumed that the parameter $\epsilon$ is chosen according to the size of the parameter $N$. If a statement is made for ``sufficiently small $\epsilon$," it should be understood that $N$ is chosen according to the size of $\epsilon$.

If an expression of the form $O(\epsilon)$ appears in an exponent, it should be assumed that the coefficient on $\epsilon$ is positive. More specifically, if we write $f(x)=g(x)^{a+O(\epsilon)}$, it means $g(x)^{a} \leq f(x)\leq g(x)^{a+C\epsilon}$, for some $C>0$. 
Similarly, for negative multiples of $\epsilon$, we will write $- O(\epsilon)$. If we write $f(x)=g(x)^{a \pm O(\epsilon)}$, it means $g(x)^{a-C\epsilon} \leq f(x)\leq g(x)^{a+C\epsilon}$, for some $C>0$.
\section{Outline}
\setcounter{subsection}{1}
The main part of the proof is carried out in Sections \ref{ContFrac}--\ref{sec:endOfArg}, eliding some technicalities which will be dealt with in later sections and the Appendix. In Section \ref{ContFrac}, we encode particular exact-order numbers via their continued fraction expansion. In Section \ref{sec:oscArg}, we indicate how particular oscillatory integrals appear in the analysis and summarise Kaufman's approach to bounding them. In Section \ref{sec:geomArg}, we introduce the core ideas of the measure construction. Section \ref{sec:combEst} shows how we can combine geometric and oscillatory data, stating a key estimate of Kaufman's. In Section \ref{sec:endOfArg}, we use these results to complete the proof of Theorem \ref{mainthm}, by making a suitable choice of scales for measure decomposition.

Section \ref{sec:CFfacts} records some basic lemmas about continued fractions we require use of throughout, in particular lemmas that concern the growth of continuants of continued fractions with suitably bounded entries. 

Section \ref{Measure} contains the technical details of the measure construction, which satisfies some suitable condition on the uniform growth of continuants. With reference to the encoding of our exact order sets, we carry out our measure construction by defining its value on admissible sequences.

Sections \ref{sec:measureOfBalls} and \ref{sec:relMeasBalls} concern geometric information about related push-forward measures, supported on the exact order numbers. Specifically, we establish the ball conditions which are used in Section \ref{sec:combEst} to complete the proof.

In order to complete Kaufman's argument, we require an approximation of certain conditional measures and we must estimate the error in approximation. Section \ref{sec:measAppr} establishes the remaining error term utilised in Section \ref{sec:combEst}. 

Appendix \ref{appendixKaufmanEsts} reproduces some key estimates of Kaufman and is included for completeness.

\section{Encoding exact-order numbers}\label{ContFrac}
\setcounter{subsection}{1}
Kaufman's argument worked via the encoding badly approximable numbers by their continued fraction expansions. The badly approximable numbers are those whose continued fraction expansion have entries bounded above by some natural number $N$. In our setting, we sparsely insert a few exceptional large partial quotients between partial quotients bounded above by some suitable large $N$. These exceptional partial quotients ensure that we have infinitely many convergents which approximate our numbers $x$ to an exact approximation order. Their sparseness will allow us to adapt the measure construction of Kaufman and Queffélec--Ramaré and their associated analysis.

\begin{mylem}\label{lem:encExact}Let $N, p \in \mathbb{N}$ and let $\mathcal{J} := \{j_k\}$ be an increasing sequence of integers Let $H := \{\eta_k\}$ be a sequence of positive real numbers that decreases to zero; for example, $\eta_k = 2^{-k}$. Define $E_{N,p,\mathcal{J},H}$ to be the set of irrational numbers $x$ having an infinite continued fraction expansion of the form 
\[[a_0; a_1, \ldots, a_{j_1 p - 1}, b_1, a_{j_1 p}, \ldots, a_{j_2, p - 1}, b_2, a_{j_2 p}, \ldots]\]
where each $a_n$ satisfies $1 \leq a_n \leq N$, and each $b_k$ satisfies 
\[b_k\in ((1 +  \eta_k/100)\rho(K_k),(1 + \eta_k/50)\rho(K_k)),\]
with $\rho$ as in \eqref{rhodef} and 
\[K_k=K(a_0,a_1,\ldots ,a_{j_1 p -1},b_1,a_{j_1 p},\ldots, a_{j_2 p-1},b_2,a_{j_2 p},\ldots ,a_{j_{k}p -1}).\] 
Then, provided that $p$ is chosen sufficiently large and the sequence $\mathcal{J}$ grows sufficiently rapidly depending on $N, p,$ and $H$, we have that $E_{N,p,\mathcal{J},H} \subset \Exact(\psi)$.
\end{mylem}
Elements of $E_{N,p,\mathcal{J},H}$ are encoded by sequences predominantly consisting of $a_n$ bounded above by $N$, but these sequences contain a sparse subsequence of entries $b_k$ that are significantly larger. These sequences (including finite initial segments) are referred to throughout as admissible sequences. The bounded partial quotients $a_n$ will be called \textbf{typical partial quotients} and will always be denoted with the letter $a$; the larger partial quotients will be called \textbf{exceptional partial quotients} and will always be denoted with the letter $b$. The letter $c$ will be used for partial quotients that may be large or small.

The next two claims will prove Lemma \ref{lem:encExact}: Claim 1 will check that there are infinitely many fractions $\frac{p_n}{q_n}$ satisfying $\left|x - \frac{p_n}{q_n} \right| \leq \psi(q_n)$, and Claim 2 will show for any $\eta > 0$ that $\left|x - \frac{p}{q} \right| \geq (1 - \eta) \psi(q)$ for all rational numbers $q$ with large $q$.
\begin{myclm}
Let $\eta > 0$. Then there exists a $Q_1(\eta)$ with the following property. Let $\frac{p_{n-1}}{q_{n-1}} = [c_0; c_1 \ldots, c_{n-1}]$ and $\frac{p_n}{q_n} = [c_0; c_1, \ldots, c_n]$ and suppose that $b$ is an integer with 
\[\left(1 + \frac{\eta}{100} \right) \rho(q_n) < b < \left(1 + \frac{\eta}{50} \right)\rho(q_n).\]
Then if $q_n > Q_1(\eta)$ and $x\in \cyl(c_0,c_1, \ldots, c_n, b)$, then 
\begin{equation}\label{pnqnbound} \left( 1 - \frac{\eta}{10} \right) \psi(q_n) \leq \left|x - \frac{p_n}{q_n}\right| \leq \psi(q_n).
\end{equation}
\end{myclm}
\begin{proof}
The continuant $q_{n+1} := K(a_0, a_1, \ldots, a_n, b)$ is given by $q_{n-1} + b q_n$. Therefore, we certainly have that 
\[q_{n-1} + b q_n \geq b q_n \geq \left(1 + \frac{\eta}{100} \right) q_n \rho(q_n),\]
and on the other hand,
\begin{IEEEeqnarray*}{Cl}
& q_{n-1} + b q_n \\
\leq & (b + 1) q_n \\
\leq & \left( \left(1 + \frac{\eta}{50} \right) \rho(q_n) + 1 \right)q_n \\
\leq & \left(1 + \frac{\eta}{40} \right) q_n \rho(q_n),
\end{IEEEeqnarray*}
if $q_n$ is sufficiently large, because $\rho(q_n) \to \infty$ by \eqref{rholimit}.
So in any case, we have
\begin{equation}\label{qnplus1bound}
\left(1 + \frac{\eta}{100} \right) q_n \rho(q_n) \leq  q_{n + 1} \leq \left(1 + \frac{\eta}{40} \right) q_n \rho(q_n)
\end{equation}
Writing $[a_0; a_1, \ldots, a_n, b]=p_{n+1}/q_{n+1}$ in reduced form and using the basic properties of continued fractions, we observe from \eqref{rhodef} that
\begin{IEEEeqnarray*}{rCl}
\left|\frac{p_n}{q_n} - x \right| & \leq & \left| \frac{p_n}{q_n} - \frac{p_{n+1}}{q_{n+1}} \right| \\
& = & \frac{1}{q_n q_{n+1}} \\
& \leq & \frac{1}{q_n \cdot \left(1 + \frac{\eta}{100} \right) q_n \rho(q_n)} \\
& \leq & \frac{1}{\left(1 + \frac{\eta}{100} \right) q_n^2 \rho(q_n)} \\
& \leq & \psi(q_n).
\end{IEEEeqnarray*}
On the other hand, writing $\frac{p_{n+2}}{q_{n+2}}$ for the convergent of $x$ occurring after $\frac{p_{n+1}}{q_{n+1}}$, the properties of continued fractions and the triangle inequality also give that
\begin{IEEEeqnarray*}{rCl}
\left|\frac{p_n}{q_n} - x \right| & \geq & \left| \frac{p_n}{q_n} - \frac{p_{n+2}}{q_{n+2}} \right| \\
& \geq & \left| \frac{p_n}{q_n} - \frac{p_{n+1}}{q_{n+1}} \right| - \left| \frac{p_{n+1}}{q_{n+1}} - \frac{p_{n+2}}{q_{n+2}} \right| \\
& \geq & \left(1 - \frac{\eta}{25} \right) \psi(q_n) - \frac{1}{q_{n+1}q_{n+2}} \\
& \geq & \left(1 - \frac{\eta}{25} \right) \psi(q_n) - \frac{1}{q_{n+1}^2} \\
& \geq & \left(1 - \frac{\eta}{25} \right) \psi(q_n) - \frac{1}{\left(1 + \frac{\eta}{100} \right)^2 (q_n \rho(q_n))^2} \\
& \geq & \left(1 - \frac{\eta}{25} \right) \psi(q_n) - \frac{\psi(q_n)}{\rho(q_n)} \\
& \geq & \left(1 - \frac{\eta}{10} \right) \psi(q_n),
\end{IEEEeqnarray*}
since $\rho(q_n) \to \infty$.
\end{proof}

\begin{myclm}
Let $\frac{p_{n-1}}{q_{n-1}} = [c_0; c_1, \ldots, c_{n-1}]$ and $\frac{p_n}{q_n} = [c_0; c_1, \ldots, c_n]$. Let $1 \leq a \leq N$ and suppose $x\in \cyl(c_0,c_1, \ldots, c_n, a)$. Then $x$ satisfies the bound
\[\frac{1}{(N + 2)} q_n^{-2} < \left|x  - \frac{p_n}{q_n} \right|.\]
In particular, we have from \eqref{q2psiqlimit} that if $q_n$ is sufficiently large depending on $N$ and $\psi$, we have
\[\psi(q_n) \leq \left|x - \frac{p_n}{q_n}\right|.\]
\end{myclm}
\begin{proof}The basic properties of the continued fraction expansion imply that $x$ satisfies the bound
\[\left| x - \frac{p_n}{q_n} \right| \geq \left| \frac{p_n}{q_n} - \frac{p_n + p_{n + 1}}{q_n + q_{n+1}}\right|;\]
therefore,
\begin{IEEEeqnarray*}{rCl}
\left| x - \frac{p_n}{q_n} \right| & \geq & \left| \frac{p_n}{q_n} - \frac{p_n + p_{n+1}}{q_n + q_{n+1}} \right| \\
& = & \left| \frac{p_n q_n + p_n q_{n+1} - p_n q_n - p_{n+1} q_n}{q_n (q_n + q_{n+1})} \right| \\
& = & \left| \frac{p_n q_{n+1} - p_{n+1} q_n}{q_n(q_n + q_{n+1})} \right| \\
& = & \frac{1}{q_n (q_n + q_{n+1})} \\
& = & \frac{1}{q_n (q_n + q_{n-1} + a q_n)} \\
& \geq & \frac{1}{(N + 2) q_n^2}
\end{IEEEeqnarray*}
as desired.
\end{proof}

\section{The oscillation argument}\label{sec:oscArg}
\setcounter{subsection}{1}
We here outline a central thread of the argument. We are interested in estimating the decay of the Fourier transform of a measure $\lambda^{\sharp}$ supported on the set of exact order numbers. In order to do so, we make a suitable decomposition of this measure at some coarse geometric scale. We will use geometric data about the measure at the coarse scale; the measure decomposition allows us to consider the effects of oscillation at fine scales. Specifics regarding the measure construction are provided in later sections.

We construct the measure, $\lambda^{\sharp}$, supported\footnote{Strictly speaking, the measure $\lambda$ we construct with is supported on a set of sequences whose image under the continued fraction map is supported in the exact order set. The associated push-forward measure $\lambda^{\sharp}$ is what we discuss in this section.} on the set of exact order numbers. This measure has a useful product structure and all elements of its support have continuants that grow in a uniform fashion. The product structure of this measure practically allows us form a conditional decomposition into probability measures defined by $\lambda^{\sharp}(\cyl(G\cdot F))=\lambda^{\sharp}(\cyl(G))\lambda^{\sharp}_G(\cyl(F))$. With $p(G)$ and $p'(G)$ defined to be the numerators corresponding to $K(G)$ and $K'(G)$, we can write
\[\hat{\lambda^{\sharp}}(\xi)=\sum_{G} \lambda^{\sharp}(\cyl(G))\int e\left(-\xi\dfrac{p(G)x+p'(G)}{K(G)x+K'(G)}\right)d\lambda^{\sharp}_{G}(x),\]
where the sum is taken over sequences of a suitable fixed length. In fact, it is a feature of our construction that the scale of $|\cyl(G)|\sim K(G)^{-2}$ will roughly be uniform over the sum, with $|\xi|^{\alpha-\epsilon}\leq K(G),K'(G)\leq |\xi|^{\alpha+\epsilon}$, where $\alpha$ is a parameter we specify in Section \ref{sec:endOfArg}.
 
In both the work of Kaufman and that of Queffélec--Ramaré, the resulting measure $\lambda^{\sharp}_{G}$ is uniform over the sequences, $G$, and has the same distribution as the original $\lambda^{\sharp}$. This is not true in our case and we will later require some technical machinery to handle this. Nevertheless, it is most instructive at this stage to suppose we can take $\lambda^{\sharp}_{G}=\mu^{\sharp}$ uniformly. Interchanging the order of summation and integration, we find that
\[\hat{\lambda^{\sharp}}(\xi)=\int \sum_{G} \lambda^{\sharp}(\cyl(G))e\left(-\xi\dfrac{p(G)x+p'(G)}{K(G)x+K'(G)}\right)d\mu^{\sharp}(x).\]
We are here considering the $\mu^{\sharp}$-average of an oscillating series, so we may expect to see some cancellation resulting from the uniformity of $\mu^{\sharp}$. The integral is not immediately amenable to treatment with the most classical tools, since we are integrating with respect to the measure $\mu^{\sharp}$.
Nevertheless, with
\[f_{\xi}(x)=\sum_{G} \lambda^{\sharp}(\cyl(G)) e\left(-\xi\dfrac{p(G)x+p'(G)}{K(G)x+K'(G)}\right),\]
the approach of Kaufman allows us to work with bounds for
\[\|f_{\xi}\|_{L^2([1,N+1])}^2=\sum_{(G,\tilde{G})} \lambda^{\sharp}(\cyl(G))\lambda^{\sharp}(\cyl(\tilde{G}))\int_1^{N+1}  e\left(-\Phi_{(G,\tilde{G})}(x)\right)dx,\]
where the phase $\Phi_{(G,\tilde{G})}$ is defined implicitly, together with a suitable ball condition on the measure $\mu^{\sharp}$.

We are now in a position to use classical oscillatory integral bounds and we divide our analysis of the sequences according to whether the phase has a critical point. There are three cases to consider: (i) the phase is nonstationary, (ii) the phase has a unique critical point, (iii) the phase is constant. The first case is where we obtain the strongest estimates. For estimate (ii), the critical point of the phase is nondegenerate and we can apply a van der Corput estimate. The size of $\xi$ and the partial quotient $K(G)$ are what determines the resulting estimate for (ii). For estimate (iii), since there is no oscillation in each of the integrals, we use geometric information about the corresponding set. More specifically, to bound (iii), we use a ball condition at the coarse scale at which we decomposed our original measure $\lambda^{\sharp}$, at the width of the intervals $\cyl(G)$, and sum the resulting estimates.

\section{The geometric argument}\label{sec:geomArg}
\setcounter{subsection}{1}
In this section, we draw attention to some of the key geometric features of the measure construction. We first carry out the construction on the space of sequences of natural numbers. The corresponding push-forward measures on the real line were those we discussed in the previous section. One important aspect of the construction is how we are able to work at uniform scales. Indeed, we construct our measure $\lambda$ to ensure that for infinite sequences $G=(c_0,c_1,\ldots)$ and $\tilde{G}=(\tilde{c}_0,\tilde{c}_1,\ldots)$ in the measure's support, we have that $K(G_n)^{-2}\approx K(\tilde{G}_n)^{-2}$, where $G_n=(c_0,c_1,\ldots,c_n)$ and $\tilde{G}_n=(\tilde{c}_0,\tilde{c}_1,\ldots,\tilde{c}_n)$.

We construct our sequences corresponding to the exact order set via the encoding discussed in Section \ref{ContFrac}. A significant difference in our work to that of Kaufman is in the introduction of exceptional partial quotients. For now, we note that at corresponding scales these ensure any probability measure can only satisfy a ball condition with a \emph{worse} exponent than the corresponding set of badly approximable numbers. To compensate for this, we include two countermeasures in the construction. Firstly, we choose our $N$, which bounds the typical partial quotients, sufficiently large. If we had no exceptional partial quotients, this choice of large $N$ would ensure the Hausdorff dimension of the badly approximable set was close to $1$ and Kaufman's original measure (which, at some scales, will behave similarly to ours) satisfies a ball condition with an exponent close to $1$. Secondly, we choose the exceptional partial quotients to be suitably sparse so that the bounded partial quotients appear only after the ball condition is near $1$. 

By construction, the measure $\lambda$ will have the following properties:
\begin{enumerate}[(i)]
\item There is an integer $p$ so that every point in the closed support of $\lambda$ has an infinite continued fraction expansion of the form
\[\bigl[(a_0; a_1, \ldots, a_{p-1}), (a_p, \ldots, a_{2p - 1}), \ldots, (a_{(j_1 - 1)p}, \ldots, a_{j_1 p - 1}), b_1,\] 
\[(a_{j_1 p}, \ldots, a_{(j_1 + 1) p - 1}), \ldots, (a_{(j_2 - 1)p}, \ldots, a_{j_2 p - 1}), b_{2}, \ldots\bigr].\]
Here, each $a_i$ is an integer from $1$ up to some fixed number $N$, and each $p$-block $\mathbf{a}_i := (a_{ip}, \ldots, a_{ip - 1})$ is independent of the previous $p$-blocks and of the $b_k$ for all $k$ such that $j_k \leq i$. 

Given the finite sequence $G_k := (\mathbf{a}_0, \ldots, \mathbf{a}_{j_1 - 1}, b_1, \mathbf{a}_{j_1}, \ldots, \mathbf{a}_{j_k - 1})$,
the exceptional entry $b_k$ is obtained as follows. We let $b_k$ be uniformly distributed in the set of integers within a fine exponential interval $I_k$, which contains only numbers approximately of scale $\rho(K(G))$. The interval $I_k$ is properly defined in Section \ref{Measure}. 

\item There is a number $\sigma > 0$  such that if $j_k \leq j < j_{k+1}$, then we have the estimate
\begin{equation}\label{typicalcontinuant}
\left| \log K(\mathbf{a}_0, \ldots, \mathbf{a}_{j_k - 1}, b_k, \mathbf{a}_{j_k}, \cdots, \mathbf{a}_{j}) - (j + (\tau - 2) j_k) \sigma \right| < \epsilon j \sigma. 
\end{equation}
\end{enumerate}
As such, we see that all sequences in the support of $\lambda$ encode exact order numbers as characterised by Lemma \ref{lem:encExact}.

Let us now observe some of the differences with Kaufman's construction. Firstly, continuant scales which occur in the construction of the measure may no longer intersect all exponential intervals, $[e^{j\sigma}, e^{(j+1)\sigma})$. As a result, we are restricted in the scales (corresponding to the geometric scale $\sim |\cyl(G)|$ in the last section) at which we can decompose our measure. Secondly, the measure we construct no longer has the self-similarity or independence properties utilised in the previous work. To deal with this, we will make a classification of the conditional measures, $\lambda_G$, and find a uniform measure of comparison over each class. 

So far, we have seen one particular scale in our analysis. This is the coarse geometric scale ($\sim |\cyl(G)|$) at which we decomposed the measure for the analysis of its Fourier transform. At this scale, we use geometric information about the measure in the application of a ball condition for the measure. Another scale important to our approach is a fine scale that allows us to relate our classical oscillatory integral estimate to our $\mu^{\sharp}$-average of an oscillating sum. The measures $\lambda_G^{\sharp}(=\mu^{\sharp})$ are probability measures defined by the relation $\lambda^{\sharp}(\cyl(G\cdot F))=\lambda^{\sharp}(\cyl(G))\lambda^{\sharp}_G(\cyl(F))$. The fine scale thus appears in two related senses, there is an \emph{absolute} fine scale, i.e. the scale of the cylinders $\cyl(G\cdot F)$, corresponding with our original measure $\lambda^{\sharp}$, and a \emph{relative} fine scale, i.e. the scale of the cylinders $\cyl(F)$, corresponding with the measures $\lambda_G^{\sharp}$. We will require a ball condition for the measures $\lambda_G^{\sharp}$ for cylinders $\cyl(F)$ at the relative fine scale, which we establish in Section \ref{sec:relMeasBalls}.

\section{Combining measure geometry and oscillation}\label{sec:combEst}
\setcounter{subsection}{1}
Recall from Section \ref{sec:oscArg} the heuristic identity
\[\hat{\lambda^{\sharp}}(\xi)=\int \sum_{G} \lambda^{\sharp}(\cyl(G))e\left(-\xi\dfrac{p(G)x+p'(G)}{K(G)x+K'(G)}\right)d\mu^{\sharp}(x).\]
For reasons we outline in the next section, we more precisely consider expressions of the form 
\[\int  \sum_{G\in \mathcal{A}_{M_1,M_2}}\lambda^{\sharp}(\cyl(G))e\left(-\xi\dfrac{p(G)x+p'(G)}{K(G)x+K'(G)}\right)d\lambda_{M_1,M_2}^{\sharp}(x),\]
with $\lambda_{M_1,M_2}=\lambda_{G_{M_1,M_2}}$ for some class representative $G_{M_1,M_2}\in \mathcal{A}_{M_1,M_2}$.

The key lemma we use for combining oscillatory data with geometric data is the following, which is the Queffélec--Ramaré version of Kaufman's lemma. 

\begin{mylem}\label{QRlem}
Let $F$ be any $C^1$ function on a compact interval $[a,b]$ satisfying $|F(t)| \leq 1$ and $|F'(t)| \leq M$. Set $m_2=\|F\|^2_{L^2([a,b])}$. Let $r > 0$ be a parameter, $0 < \beta < 1$ be an exponent, and suppose that $\mu$ is any probability measure on $[a,b]$ that satisfies the property that for every ball $B$ of radius $r/M$, $\mu$ satisfies the $\beta$-ball condition
\[\mu(B) \leq |B|^{\beta}.\]
Then we have the estimate
\[\int_a^b |F(x)| \, d \mu(x) \leq 2 r + (r/M)^{\beta} (1 + m_2 M r^{-3}).\]
\end{mylem}
The proof of Lemma \ref{QRlem} is presented alongside the proofs of the necessary input data, Lemmas \ref{Mlemma} and \ref{m2lemma} in Appendix \ref{appendixKaufmanEsts}. They are included for completeness but are essentially the same as in those previous works of Kaufman and Queffélec--Ramaré. Briefly, Lemma \ref{QRlem} is proven by a pigeonhole argument. We can decompose $[a,b]$ into intervals, $I$, of width $r/M$ and the number of intervals $I$ for which $\inf_{x\in I}|F(x)| \geq r$ is controlled using $m_2$. By the control on $F'$, the remaining intervals $I'$ must have $\sup_{x\in I}|F(x)|\leq 2r$.

In our case, we consider
\begin{equation}\label{Fdef}F(x)=f_{\xi}(x)=\sum_{G\in \mathcal{A}_{M_1,M_2}} \lambda^\sharp(\cyl(G)) e\left(-\xi\dfrac{p(G)x+p'(G)}{K(G)x+K'(G)}\right)\end{equation}
and $\|f_{\xi}\|_{L^2([1,N+1])}^2$ can be treated as a sum of classical oscillatory integrals. The sequences $G$ over which we decompose our measure will be chosen such that $|\xi|^{\alpha-\epsilon}\leq K(G) \leq |\xi|^{\alpha+\epsilon}$, 
for some suitable parameter $\alpha$ to be considered in more detail in Section \ref{sec:endOfArg}.

\begin{mylem}\label{Mlemma}
Let $F$ be as in \eqref{Fdef}. Then provided that $|\xi|$ is large enough, we have the estimate
\[M := \max_{x \in [1, N + 1]} |F'(x)| \leq |\xi|^{1 - 2 \alpha + 3 \epsilon}.\]
\end{mylem}

The oscillatory integral argument, which was outlined in Section \ref{sec:oscArg}, hews closely to the work of Kaufman and allows us to establish the following lemma.
\begin{mylem}\label{m2lemma}
Let $F$ be as in \eqref{Fdef}. Then provided that $|\xi|$ is large enough, we have the estimate
\[m_2 = \int_1^{N+1} |F(x)|^2 \, dx \leq |\xi|^{ \frac{3\alpha-1}{2}  + O( \epsilon)} + |\xi|^{-\frac{\tau \alpha}{\tau - 1} +
 O( \epsilon)}.\]
\end{mylem}

\section{Completing the argument}\label{sec:endOfArg}
\setcounter{subsection}{1}
The work of Queffélec--Ramaré makes clear one of the tolerances in Kaufman's argument is in the choice of coarse and fine scales. By making an optimal choice adapted to $\dim\textbf{Bad}'(S)$, (as well as refinements to the stationary phase estimates), they improve on Kaufman's result. For our purposes, this insight is useful as it suggests we make an adaptive choice of scales so that (i) we can actually decompose our measure at the desired coarse scale and (ii) we have a suitable ball condition at both of the coarse and fine scales to run our analysis. Nevertheless, we do not optimise the choice of scales. 

When describing scales, we will use the term \textbf{continuant scale} to refer to the approximate size of the continuant of a finite continued fraction, and the term \textbf{geometric scale} to refer to the approximate Lebesgue measure of an interval. A continuant scale $\zeta$ corresponds to a geometric scale of $\zeta^{-2}$. We now consider which are those typical scales our measure enables us to work at and how to address the exceptional scales.

\begin{mydef}\label{def:typExcScales}Let $\zeta \gg 1$ be a real number. We will say that the (continuant) scale $\zeta$ is \textit{exceptional} if there exists some $k \in \mathbb{N}$ such that $(1 - 2 \epsilon) j_k \sigma < \log \zeta < (\tau - 1 + 2 \epsilon) j_k \sigma$. Otherwise, we say that the scale $\zeta$ is \textit{typical}. With $(\mathbf{a}_0, \mathbf{a}_1, \ldots) \in \supp \lambda$, if $k$ is such that $(\tau - 1 + 2 \epsilon) j_k \sigma < \log \zeta < (1 - 2 \epsilon) j_{k+1} \sigma$, then we define
\begin{equation}\label{eq:typjScale}j(\zeta) :=\left \lfloor  \dfrac{\log\zeta - (\tau - 2) j_k\sigma}{\sigma} \right\rfloor.\end{equation}\end{mydef}
The point is that if a scale is typical, then every element in the support of $\lambda^\sharp$ has a partial quotient with continuant close to $\zeta$ (see Lemma \ref{lem:scales}). 

We will mimic Kaufman's strategy to estimate $\widehat{\lambda^{\sharp}}(\xi)$ by considering a coarse continuant scale $\zeta=|\xi|^{\alpha}$ at which to decompose our measure, as discussed. In the work of Kaufman and Queffélec--Ramaré, the measure decomposition is carried at continuant scale $\zeta=|\xi|^{\alpha}$, for some suitable $\alpha$. Kaufman considers the scale $\alpha=0.17$ uniformly across the class of measures and Queffélec--Ramaré make a choice of $\alpha$ which is sensitive to the dimension of the measure. In our case, for any fixed $\alpha$, $\zeta = |\xi|^{\alpha}$ may be exceptional. As such, we decompose our measure at scale $|\xi|^{\alpha}$ with variable $\alpha$ chosen to ensure the scale is typical. In fact, we only require two choices $\alpha_0$ and $\alpha_1$. If $|\xi|^{\alpha_0}$ is an exceptional scale, our choice of $\alpha_1$ will ensure that $|\xi|^{\alpha_1}$ is typical and we can run the argument by decomposing our measure at this continuant scale. In compensation for decomposing at the larger scale $\alpha_1$, we are able to utilise a stronger relative ball condition at the fine scale, with exponent closer to $1$, rather than a weaker scale-uniform ball condition.\footnote{The measures considered by Kaufman and Queffélec--Ramaré are more aptly treated with a scale-uniform ball condition, however Kaufman uses a uniform ball condition over all scales \emph{and} measures and this is one point by which Queffélec--Ramaré are able to obtain improved estimates.}

For our purposes, it is suitable to choose $\alpha_0 := \frac{10 - \sqrt{73}}{9} \approx 0.162$. Furthermore, we define the exponent $\alpha_1$ by the equation
\begin{equation}\label{alphaprimedef}
\alpha_1 = (\tau - 1 + 10 \epsilon) \alpha_0.
\end{equation}
It is easy to see for $\tau < (13+\sqrt{73})/8$ that $0.161 \leq \alpha_1 \leq 0.274$. As we will see, this choice of $\alpha_1$ ensures that if $|\xi|^{\alpha_0}$ is an exceptional scale, then $|\xi|^{\alpha_1}$ is typical.

\begin{mylem}\label{lem:scales}Let $\xi\in\mathbb{R}$ be large. For some $0.161\leq\alpha\leq 0.274$, suppose that $\zeta=|\xi|^{\alpha}$ is a typical scale. Then, for an admissible sequence $G$ containing $j(\zeta)$ $\mathbf{a}_i$'s, we have that
\[K(G)\in [|\xi|^{\alpha-\epsilon}, |\xi|^{\alpha+\epsilon}]\quad\text{and}\quad|\cyl(G)|\in \left[|\xi|^{-2\alpha-2\epsilon},|\xi|^{-2\alpha+2\epsilon}\right].\] \end{mylem}
\begin{proof}
Suppose that $j_k\leq j(\zeta)<j_{k+1}$. We then have by an application of \eqref{typicalcontinuant}, followed by the assumption $ j_k >\sigma\epsilon^{-1}$ and a further application of \eqref{typicalcontinuant} that
\[\left|\log\zeta- \log K(G)\right|\]
\[\leq 
\left|\log\zeta- (j(\zeta) +j_k (\tau - 2))\sigma\right|+\left|\log K(\mathbf{a}_0, \mathbf{a}_1, \ldots,\mathbf{a}_{j(\zeta)})- (j(\zeta) +j_k(\tau - 2)) \sigma\right|\] 
\[\leq \sigma+\epsilon j(\zeta)\sigma<2\epsilon\sigma j(\zeta)\leq \frac{2\epsilon j(\zeta)\log K(G)}{j(\zeta) +j_k(\tau - 2)-\epsilon j(\zeta) }\leq 3\epsilon \log K(G).\] 
As such, we have that \[\zeta\in \left[K(G)^{1-3\epsilon}, K(G)^{1+3\epsilon}\right]\] and \[K(G)\in \left[\zeta^{\frac{1}{1+3\epsilon}}, \zeta^{\frac{1}{1-3\epsilon}}\right]\subset \left[\zeta^{1-3.1\epsilon}, \zeta^{1+3.1\epsilon}\right].\]
We can also observe that since $|\cyl(G)|\sim K(G)^{-2}$, we have that 
\[|\cyl(G)|\in \left[\zeta^{-2-7\epsilon}, \zeta^{-2+7\epsilon}\right],\]
where we have absorbed the constant of comparison with an additional $\zeta^{0.8\epsilon}$ loss. We choose $\alpha\in[0.161,0.274]$, so the result follows. 
\end{proof}

Note that if $K(G)$ is sufficiently large, for an admissible sequence $G$, we also have
\begin{equation}\label{denominatorratio}
1 + \frac{1}{N + 2} \leq \frac{K(G)}{K'(G)} \leq N + 1.
\end{equation}
We require a use of the following classification of measures in order to interchange the order of integration and summation in the manner outlined in Section \ref{sec:oscArg}.

\begin{mydef}\label{def:AM1M2class}
We form a minimally sized cover of the rectangle $[|\xi|^{\alpha - \epsilon}, |\xi|^{\alpha + \epsilon}]\times [|\xi|^{\alpha - \epsilon}, |\xi|^{\alpha + \epsilon}]$ by disjoint squares of length $|\xi|^{\alpha - 200 \epsilon}$. There are $\lesssim |\xi|^{600 \epsilon}$ such boxes.

Suppose $(M_1, M_2)$ is the bottom-left corner of one of these boxes. Say $G \in \mathcal{A}_{M_1, M_2}$ if the point $(K'(G), K(G))$ lies inside the box with bottom-left corner $(M_1, M_2)$; i.e. $\mathcal{A}_{M_1,M_2}$ consists of those values of $G$ satisfying the pair of inequalities
\begin{IEEEeqnarray*}{rCcCl}
M_1 & \leq & K(G) & \leq & M_1 + |\xi|^{\alpha - 200 \epsilon} \\
M_2 & \leq & K'(G) & \leq & M_2 + |\xi|^{\alpha - 200 \epsilon}.
\end{IEEEeqnarray*}
If such a $G$ exists, then, for large $\xi$, equation \eqref{denominatorratio} implies that $M_1$ and $M_2$ must satisfy
\begin{equation}\label{M1M2ratio}
1 < \frac{M_1}{M_2} \leq N + 1.1.
\end{equation}
For each pair $(M_1, M_2)$ for which $\mathcal{A}_{M_1, M_2}$ is nonempty, choose a representative element $G_{M_1, M_2}$. Write $\lambda_{M_1, M_2}$ for $\lambda_{G_{M_1, M_2}}$.
\end{mydef}

\begin{proof}[Proof of Theorem \ref{mainthm}]
Suppose that $\zeta=|\xi|^\alpha$ is a typical scale and let $j(\zeta)$ be as in \eqref{eq:typjScale}. Let $\mathcal{S}(\zeta)$ denote the family of admissible sequences $G$ containing $j(\zeta)$ $\mathbf{a}_i$'s. Consider $z\in\cyl(G)\cap \supp\lambda^{\sharp}$. We know that \[z=\frac{p(G)x+p'(G)}{q(G)x+q'(G)}=[\mathbf{a}_0,\ldots, \mathbf{a}_{j(\zeta)},x].\]
That is, $x$ is the real number with the continued fraction expansion which is the tail of the expansion of $z$. Given $z\in \cyl(G)$, such numbers $x$ are distributed according to the measure $\lambda_G^{\sharp}$ (see Section \ref{sec:relMeasBalls} and, in particular, \eqref{contfracshift}). By making the associated change-of-variables, we can write 
\begin{equation}\label{measuredecomp1}
\widehat{\lambda^{\sharp}}(\xi) = \sum_{G \in \mathcal{S}(\zeta)} \lambda^{\sharp}(\cyl(G)) \int e \left(\xi  \frac{p(G)x + p'(G)}{q(G)x + q'(G)} \right) \, d \lambda_G^{\sharp}(x).
\end{equation}
Observe that the coefficients $\lambda^{\sharp}(\cyl(G))$ sum to $1$.

By a classification of the the conditional measures and comparison, the precise identity we make use of has the form
\begin{IEEEeqnarray*}{rCl}\hat{\lambda^{\sharp}}(\xi)&=&\sum_{(M_1,M_2)} \int  \sum_{G\in \mathcal{A}_{M_1,M_2}}\lambda^{\sharp}(\cyl(G))e\left(-\xi\dfrac{p(G)x+p'(G)}{K(G)x+K'(G)}\right)d\lambda_{M_1,M_2}^{\sharp}(x)\\
+&&\sum_{(M_1,M_2)} \int  \sum_{G\in \mathcal{A}_{M_1,M_2}}\lambda^{\sharp}(\cyl(G))e\left(-\xi\dfrac{p(G)x+p'(G)}{K(G)x+K'(G)}\right)d\left(\lambda_G^{\sharp}-\lambda_{M_1,M_2}^{\sharp}\right)(x).\IEEEyesnumber \label{eq:measDecompMain}
\end{IEEEeqnarray*}

The first term in \eqref{eq:measDecompMain} is treated as outlined in the previous sections, combining oscillation estimates with a ball condition for $\lambda_{M_1,M_2}$ at a fine scale. The additional sum over $(M_1,M_2)$ has $\sim |\xi|^{600 \epsilon}$ many terms, which our decay estimates will compensate for. 
Lemma \ref{classlemma}, characterising the measure classification error, tells us that the second term in \eqref{eq:measDecompMain} is $O(|\xi|^{-\epsilon})$. The proof is completed by the combination Lemma \ref{classlemma}, the count $|\mathcal{A}_{M_1,M_2}|\sim |\xi|^{600 \epsilon}$ from Definition \ref{def:AM1M2class}, and the below Claim \ref{pieceestimate} for a sufficiently small choice of $\epsilon$ to ensure $-c_\tau+O(\epsilon)<-700\epsilon$. 

\begin{myclm}\label{pieceestimate}
Suppose that $\xi$ is sufficiently large. In the case that $|\xi|^{\alpha_0}$ is a typical scale, set $\alpha=\alpha_0$. Otherwise, set $\alpha=\alpha_1$. Then, in either case, $|\xi|^{\alpha}$ is a typical scale. Furthermore, for admissible sequences $G$ containing $j(|\xi|^{\alpha})$ $\mathbf{a}_i$'s, we have that, for each $M_1$ and $M_2$,
\[\int \sum_{G \in \mathcal{A}_{M_1, M_2}} \lambda^\sharp(\cyl (G)) e \left( \xi  \frac{p(G)x + p'(G)}{q(G)x + q'(G)} \right) \, d \lambda_{M_1, M_2}^\sharp(x) \lesssim  |\xi|^{-700\epsilon} + |\xi|^{-c_\tau + O(\epsilon)},\]
for some $c_{\tau}>0$.
\end{myclm}
\begin{proof}[Proof of Claim \ref{pieceestimate}]
First, note that the measures $\lambda_{M_1,M_2}=\lambda_{G_{M_1,M_2}}$ for some $G_{M_1,M_2}\in\mathcal{A}_{M_1,M_2}$ satisfy a ball condition at all scales with exponent $0.4\leq \beta_F\leq 1$. This is a consequence of Lemma \ref{lem:relCylBad}, which gives \[\beta_F =1-\frac{1-2 / \tau }{1 - 2 \alpha}-O(\epsilon).\]
Indeed, we know that $2\leq \tau < \frac{13 + \sqrt{73}}{8}=\bar{\tau}$ and $0.161\leq\alpha\leq 0.274$ so that 
\[\beta_F\geq 1-(1-2 / \bar{\tau} )/(1 - 2 \cdot 0.274)-0.01\geq 0.4.\]

We consider specifically what is given when we insert Lemmas \ref{Mlemma} and \ref{m2lemma} into Lemma \ref{QRlem}. We make a choice for $r$ which is far from optimal, since we only require \emph{some} decay in $|\xi|$. If we choose $r=|\xi|^{- 1000\epsilon}$, the geometric scale $r/M$ turns out to be $ |\xi|^{-1+2\alpha-O(\epsilon)}$. Let $r =|\xi|^{-1000\epsilon}$ and suppose that $\mu$ is any probability measure on $[1,N+1]$ that satisfies Lemmas \ref{Mlemma} and \ref{m2lemma} and is subject to a ball condition of exponent $1\geq \beta_F\geq 0.4$ at geometric scale $r/M= |\xi|^{-1+2\alpha - O(\epsilon)}$. Then we can obtain the estimate 
\[
\int_1^{N+1} |f_{\xi}(x)| \, d \mu(x)  \leq |\xi|^{-700\epsilon} + |\xi|^{ O(\epsilon)}M^{1-\beta_F}m_2.
\]
Indeed, the $2r$ and $(r/M)^{\beta_F}$ terms from Lemma \ref{QRlem} are together bounded by the first $ |\xi|^{-700\epsilon}$ summand. Up to a factor of $r^{\beta_F-3}=|\xi|^{O(\epsilon)}$, it remains to sum $M^{1-\beta_F}m_2$. 

\paragraph*{Case 1: $|\xi|^{\alpha_0}$ is a typical scale} Lemma \ref{lem:relCylBad} tells us that $\lambda_{M_1,M_2}^\sharp=\lambda_{G_{M_1,M_2}}^\sharp$ satisfies a ball condition with exponent $\beta_F := \frac{2/\tau - 2 \alpha }{1 - 2 \alpha} - O(\epsilon)$ at the fine geometric scale $|\xi|^{-1+2\alpha}$. Observe that if $\beta = \frac{2/\tau - 2 \alpha }{1 - 2 \alpha} - O(\epsilon)$, then $M^{1- \beta} \leq |\xi|^{1 - 2/\tau + O(\epsilon)}$. 

So
\[ m_2 M^{1 - \beta} \leq \left(|\xi|^{\frac{3\alpha-1}{2} } + |\xi|^{-\frac{\tau \alpha}{\tau - 1}} \right) |\xi|^{1 - \frac{2}{\tau} + O(\epsilon)} \]
In order for this to be smaller than $|\xi|^{-c_\tau + O(\epsilon)}$ for some $c_\tau>0$, it is sufficient to show this for each summand. That is, it is sufficient for $\alpha$ to satisfy the inequalities
\begin{IEEEeqnarray}{rCl}
\frac{3\alpha-1}{2} + \left(1 - \frac{2}{\tau} \right) & < & 0 \label{typicalineq1} \\
- \frac{\tau \alpha}{\tau - 1} + \left(1 - \frac{2}{\tau} \right) & < & 0. \label{typicalineq2}
\end{IEEEeqnarray}
Solving \eqref{typicalineq1} for $\alpha$ gives
\begin{equation}\label{typicalineq1rewrite}
\alpha < -\frac{1}{3} + \frac{4}{3 \tau}.
\end{equation}
Observe that the right side is a decreasing function of $\tau$, so if this inequality is satisfied for a given $\alpha$ and $\tau$, it is also satisfied if $\alpha$ is held constant and $\tau$ is decreased.

Solving \eqref{typicalineq2} for $\alpha$ gives
\begin{equation}\label{typicalineq2rewrite}
\alpha > \frac{\tau^2 - 3 \tau + 2}{\tau^2}   = 1 - \frac{3}{\tau} + \frac{2}{\tau^2}.
\end{equation}

The right side is a strictly increasing function of $\tau$ for $\tau \geq 2$ (in fact for $\tau \geq 4/3$), so if this inequality is satisfied for some $\alpha$ and $\tau$, it is also satisfied if $\alpha$ is held constant and $\tau$ is decreased. As such, we can verify that \eqref{typicalineq1rewrite} and \eqref{typicalineq2rewrite} hold when $\tau < \frac{13 + \sqrt{73}}{8}=\overline{\tau}$ and  $\alpha_0 = \frac{10 - \sqrt{73}}{8}$. In order to check this, it is enough to check for these values of $\overline{\tau}$ and $\alpha_0$ that
\[\alpha_0 = \frac{\overline{\tau}^2 - 3 \overline{\tau} + 2}{\overline{\tau}^2} = -\frac{1}{3} + \frac{4}{3 \overline{\tau}},\]
which is easily verified. Provided we make a choice \[0>-c_{\tau}\geq \max\left\lbrace \frac{3\alpha_0-1}{2}+(1-\frac{2}{\tau}), - \frac{\tau \alpha_0}{\tau - 1} + \left(1 - \frac{2}{\tau} \right)\right\rbrace,\]
the proof is complete. 

\paragraph*{Case 2: $|\xi|^{\alpha_0}$ is an exceptional scale} First, we will verify that $|\xi|^{\alpha_1}$ is a typical scale. Observe that for $\tau \leq 4$ we have
\begin{IEEEeqnarray*}{rCl}
\log(|\xi|^{\alpha_1}) & = & (\tau - 1 + 10 \epsilon) \log(|\xi|^{\alpha_0}) \\
& \geq & (\tau - 1 + 10 \epsilon)(1 - 2 \epsilon) j_k \sigma \\
& \geq & (\tau - 1 - 2 \epsilon - 2 (\tau - 2) \epsilon + 10 \epsilon - 20 \epsilon^2) j_k \sigma \\
& \geq & (\tau - 1 + 3 \epsilon) j_k \sigma. \IEEEyesnumber \label{xialphatypical1}
\end{IEEEeqnarray*}
On the other hand, if $j_{k+1}$ is sufficiently large depending on $j_k$, we have 
\begin{equation}\label{xialphatypical2}
\log |\xi|^{\alpha_1} \ll \frac{1}{100} j_{k+1} \sigma.
\end{equation}
In particular, \eqref{xialphatypical1} and \eqref{xialphatypical2} imply that $\zeta=|\xi|^{\alpha_1}$ is a typical scale.  Moreover, provided that $j_{k+1}$ is sufficiently large depending on $j_k$, there are no exceptional partial quotients corresponding to a continuant $K$ with $|\xi|^{\alpha_1} \leq K \leq |\xi|^{10}$. 

In this case, we can apply Lemma \ref{lem:relCylGood} to conclude that $\lambda_{M_1, M_2}^\sharp$ satisfies a ball condition with the improved exponent $\beta = 1 - O(\epsilon)$ at the relevant geometric scale $|\xi|^{1 - 2 \alpha_1 +  O(\epsilon)}$.

As in the case where $|\xi|^{\alpha_0}$ is a typical scale, it is enough to establish that $m_2 M^{1 - \beta} < |\xi|^{-c_{\tau} + O(\epsilon)}$. We have the same sort of bound as in the typical case, but with $\beta = 1 -  O(\epsilon)$. Observing that since $1 - \beta =  O(\epsilon)$, we have the estimate
\[m_2 M^{1 - \beta} \leq |\xi|^{O (\epsilon)} (|\xi|^{\frac{3 \alpha-1}{2}+ O (\epsilon)} + |\xi|^{-\frac{\tau \alpha}{\tau - 1}  + O( \epsilon)})\]
As before, it is sufficient to see that each term is bounded above by $|\xi|^{-c_{\tau} + O(\epsilon)}.$
This leads to the pair of inequalities
\begin{IEEEeqnarray}{rCl}
\frac{3 \alpha-1}{2} & < &0\label{exceptionalineq1} \\
- \frac{\tau \alpha}{\tau - 1}  & < & 0 \label{exceptionalineq2}.
\end{IEEEeqnarray}
Solving \eqref{exceptionalineq1} for $\alpha$ gives 
\begin{equation}\label{exceptionalineq1rewrite}
\alpha < 1/3.
\end{equation}

On the other hand, solving \eqref{exceptionalineq2} for $\alpha$ gives
\begin{equation}\label{exceptionalineq2rewrite}
\alpha > 0.
\end{equation}
Recall that we chose $\alpha_0 = \frac{10 - \sqrt{73}}{9}$, and that we defined $\alpha_1$ in \eqref{alphaprimedef}. We have $\alpha_1 = (\tau - 1)\alpha_0 + O(\epsilon)$. We are interested in the case in which $\tau< \frac{13 + \sqrt{73}}{8}$, we observe that with these choices that for sufficiently small $\epsilon$ we have $0.161 \leq \alpha_1 \leq 0.274$. Hence, $\alpha_1$ will satisfy \eqref{exceptionalineq1} and \eqref{exceptionalineq2}, as desired. Provided we make a choice
\[0>-c_{\tau}\geq \max\left\lbrace \frac{3(\tau - 1)\alpha_0-1}{2}, - \tau \alpha_0 \right\rbrace,\]
the proof is complete. 

Having considered Case 1 and Case 2, we see that an appropriate choice of $-c_\tau$ is given by 
\[ \max\left\lbrace \frac{3\alpha_0-1}{2}+(1-\frac{2}{\tau}), - \frac{\tau \alpha_0}{\tau - 1} + \left(1 - \frac{2}{\tau} \right),\frac{-1+3(\tau - 1)\alpha_0}{2}, - \tau \alpha_0 \right\rbrace<0.\]
\end{proof}
\end{proof}

\begin{myrmk}
The potential improvements suggested by Remark \ref{mainthmrmk} and the work of Queffélec--Ramaré would most obviously be obtained by a better choice of scales $|\xi|^{\alpha}$ in this part of the argument, sensitive to their relation with ball conditions $\mu(B)\leq |B|^{\beta}$ for fine-scale balls. However, the precise interaction of the parameters $\alpha$ and $\beta$ requires further technical calculations (for Kaufman and Queffélec--Ramaré, $\beta$ is scale-uniform for a given measure). It is also possible that improved stationary phase estimates are possible. However, our measure structure complicates a direct adaptation of Queffélec--Ramaré's improved estimate.
\end{myrmk}

\section{Continued fraction lemmas}\label{sec:CFfacts}
\setcounter{subsection}{1}
We will need some basic lemmas concerning continued fractions. The first is a trivial modification of one of the lemmas in the Kaufman paper. 

\begin{mylem}[Gluing lemma]\label{gluinglemma}
Let $K(c_0, \ldots, c_n)$ and $K(d_0, \ldots, d_m)$ be the denominators of any two finite continued fractions, and suppose that $1 \leq d_0 \leq N$. Then there is a constant $C_N$ depending only on $N$ such that
\[\left| \log K(c_0, \ldots, c_n, d_0, \ldots, d_m) - \log K(c_0, \ldots, c_n) - \log K(d_0, \ldots, d_m) \right| \leq C_N.\]
\end{mylem}
\begin{proof}
Writing $\frac{p_{n-1}}{q_{n-1}}$ and $\frac{p_n}{q_n}$ for the final two convergents of $[c_0; \ldots, c_n]$ and writing $\frac{p^*}{q^*}$ for $[d_0; \ldots, d_m]$, we have that
\begin{IEEEeqnarray*}{rCl}
[c_0; \ldots, c_n, d_0, \ldots, d_m] & = & \frac{p_{n-1} + \frac{p^*}{q^*} p_n}{q_{n-1} + \frac{p^*}{q^*} q_n}. \\
& = &  \frac{q^* p_{n-1} + p^* p_n}{q^* q_{n-1} + p^* q_n}. \label{combinedcontfrac} \IEEEyesnumber
\end{IEEEeqnarray*}
It is a well-known property of continued fractions that the fraction \eqref{combinedcontfrac} will always be reduced. So it follows that 
\[K(c_0, \ldots, c_n, d_0, \ldots, d_m) = q^* q_{n-1} + p^* q_n.\]
But we have that $q^* = K(d_0, \ldots, d_m)$ and that $q^* \leq p^* \leq (N + 1) q^*$ and $q_{n-1} \leq q_n$.  So
\[q^* q_n \leq K(c_0, \ldots, c_n, d_0, \ldots, d_m) \leq (N + 2) q^* q_n \]
Taking logarithms gives the desired inequality.
\end{proof}

We now record two simple lemmas.
\begin{mylem}For any finite sequence of natural numbers $G$, we have that $K(G)^{-2}\sim |\cyl(G)|.$
\end{mylem}
\begin{proof}We denote by $g$ the continued fraction map and set $G=(c_0,c_1,\ldots,c_n)$ and $\tilde{G}=(c_0,c_1,\ldots,c_n+1)$. We have that $|\cyl(G)|=|g(G)-g(\tilde{G})|=|g(G)-g(G\cdot 1)|=1/(K(G)(K(G)+K'(G)))\sim K(G)^{-2}$. 
\end{proof}

\begin{mylem}Consider the finite bounded sequence $\mathbf{a}=(a_1,a_2,\ldots,a_j)$ for some $j\leq p<\infty$, with $1\leq a_i\leq N$.  Then, for any finite sequence $G$, we have that $K(G\cdot \mathbf{a})\sim_{N,p} K(G)$.
\end{mylem}
\begin{proof}This is a simple consequence of the continued fraction algorithm since $K(\tilde{G}\cdot a_i) = a_iK(\tilde{G})+K'(\tilde{G})\in [K(\tilde{G}), (N+1)K(\tilde{G})]$.
\end{proof}

\section{Construction of the measure}\label{Measure}
\setcounter{subsection}{1}
We are now ready to proceed with the measure construction. Recall our use of a small parameter $\epsilon$ in the previous sections. In particular, the a priori bounds for the measures of balls used in Section \ref{sec:endOfArg} required a choice of $\epsilon = \epsilon(\tau)$ sufficiently small to obtain a decay estimate. Note that such a choice of $\epsilon$ should now be regarded as fixed, since other parameters $N$, $m$, $p$ fixed in this section depend on our choice of $\epsilon$. It remains to establish the actual ball conditions for the measures we construct, but they are derived and expressed uniformly over our choice of parameters and so there is no contradiction in regarding the parameter $\epsilon$ as being fixed at this stage.

We choose $N$ sufficiently large to ensure that $1-\epsilon<\dim\textbf{Bad}(N)<1-\epsilon/10$

Following Kaufman, we define the sum
\[S_m := \sum_{a_0, a_1, \ldots, a_{m-1} = 1}^N K(a_0, \ldots, a_{m-1})^{-2 (1 - \epsilon)},\]
and we define a probability measure $\nu_m$ on $\{1, \ldots, N\}^m$ by
\begin{equation}\label{eq:numMeasDef}\nu_m(a_0, \ldots, a_{m-1}) =  \frac{1}{S_m} K(a_0, \ldots, a_{m-1})^{-2 (1 - \epsilon)}\end{equation}
Then $\nu_m$ is a probability measure on $m$-tuples $G = (a_0, \ldots, a_{m-1})$ that captures the $(1 - \epsilon)$-dimensional Hausdorff content of $\cyl(G)$. Since we know that $\dim\textbf{Bad}(N)>1-\epsilon$, $S_m\rightarrow \infty$ as $m\rightarrow \infty$. It should be noted that we require a suitably large choice of $m$ to account for the constants arising from the application of the gluing lemma, Lemma \ref{gluinglemma}. This will be evident in what follows.

By applying the law of large numbers (the central limit theorem would be sharper but its use is not required for our purposes) we construct, according with Kaufman, a probability measure $\bar \nu_p$ on $\{1, \ldots, N\}^p$, where $p = Jm$ for some large $J$ (depending on $m$). Provided $J$ is chosen sufficiently large, we distribute the mass of $\nu\times\nu\times\ldots\times\nu$ over a suitable set $E$ with $\nu\times\nu\times\ldots\times\nu(E)\geq \frac{1}{2}$, setting $\bar\nu_p (S)=\nu\times\nu\times\ldots\times\nu(E\cap S)/\nu\times\nu\times\ldots\times\nu(E)$, so that $\bar \nu_p$ has two key properties:
\begin{enumerate}[(a)]
\item The measure $\bar \nu_p$ satisfies the bound
\[\bar \nu_p(a_0, \ldots, a_{Jm - 1}) \leq 2 \prod_{i=1}^{J} \nu_m(a_{(i-1)m}, \ldots, a_{(i-1)m+m-1})\]
\item There exists a number $\sigma > 0$ such that for all $(a_0, \ldots, a_{Jm - 1}) \in \supp \bar \nu_p$ we have
\[\left| \log K(a_0, \ldots, a_{Jm - 1}) - \sigma \right| < \frac{\epsilon}{500} \sigma.\]
\end{enumerate}

Following Kaufman, we define the measure $\lambda$. Let $\eta_k$ be a sequence of small, positive real numbers that decrease to zero; for example, $\eta_k = 2^{-k}$ will suffice. The measure $\lambda$ will be supported on integer sequences of the form 
\[(\mathbf{a}_0, \mathbf{a}_1, \ldots, \mathbf{a}_{j_1- 1}, b_{1}, \mathbf{a}_{j_1}, \ldots) \]
where we have that each $\mathbf{a}_j$ is the $p$-tuple $(a_{jp}, \ldots, a_{(j+1)p - 1})$, where $p = Jm$, $\{j_k\}_{k=1}^{\infty}$ is a suitable rapidly-growing sequence,\footnote{The sequence $j_k$ must grow sufficiently fast to ensure that the partial quotients encode exact order numbers as in Section \ref{ContFrac}. Furthermore, rapid growth of the $j_k$ is necessary to ensure the measure satisfies suitable ball conditions, as will be seen in Section \ref{sec:measureOfBalls}.} and each integer $a_l$ lies in $\{1, \ldots, N\}$ and each integer $b_k$ satisfies
\begin{equation}\label{eq:bkInt}b_k \in T_k := \mathbb{N} \cap I_{\eta_k}(G_k),\end{equation}
where $G_k := (\mathbf{a}_0, \ldots, \mathbf{a}_{j_1 - 1}, b_1, \mathbf{a}_{j_1}, \ldots, \mathbf{a}_{j_k - 1})$ and
\[I_{\eta_k}(G_k)=\left[ \left(1 + \frac{\eta_k}{1000} \right)^{\gamma_{\eta_k}(G_k)}, \left( 1 + \frac{\eta_k}{1000} \right)^{\gamma_{\eta_k}(G_k) + 1} \right)\]
is an interval contained in $\left( \left(1 + \frac{\eta_k}{100} \right) \rho(G_k), \left(1 + \frac{\eta_k}{50} \right) \rho(G_k) \right)$. More precisely, $I_{\eta_k}(G_k)$ is given as above when we define the quantity $\gamma_{\eta_k}(G_k)$ by
\begin{equation}\label{eq:gammaEtaDef}\gamma_{\eta_k}(G) := \left\lfloor \frac{\log \left(\rho(K(G)) \left(1 + \frac{\eta_k}{75} \right)\right)}{\log (1 + \frac{\eta_k}{1000})} \right\rfloor.\end{equation}

It is easy to see for this choice of $\gamma_{\eta_k}(G_k)$ that $I_{\eta_k}(G_k)$ will lie inside $\left( \left(1 + \frac{\eta_k}{100} \right) \rho(G_k), \left(1 + \frac{\eta_k}{50} \right) \rho(G_k) \right)$. We will define the space of sequences of this type by $\mathcal{S}(\psi, \{\eta_k\})$.

We define the measure $\nu_m$ as in \eqref{eq:numMeasDef}. To define the probability measure $\lambda$, we use the following procedure. Given a finite sequence $G$ of integers, recall that $G^*$ is the set of infinite sequences whose first $|G|$ terms agree with $G$. We will describe the measure $\lambda$ by describing its value on such sets $G^*$, where $G$ is an admissible sequence in $\mathcal{S}(\psi, \{\eta_k\})$. The measure is $0$ elsewhere.

We define the measure $\lambda$ as follows. 
\begin{itemize}
\item If $G = (\mathbf{a}_0, \ldots, \mathbf{a}_{j})$ for $0\leq j <j_1$, define 
\[\lambda(G^*) = \bar \nu_p(\mathbf{a}_0) \times \cdots \times \bar \nu_p(\mathbf{a}_{j}).\]
\item If $G = (\mathbf{a}_0, \cdots, \mathbf{a}_{j_k - 1}, b_k)$ for some $k$, then we take $\tilde G=(\mathbf{a}_0, \ldots, \mathbf{a}_{j_1 - 1})$ and we define 
\[\lambda(G^*) = \lambda(\tilde G)  \times |T_k|^{-1},\]
with the integer interval $T_k$ as in \eqref{eq:bkInt}.
\item If $G = (\mathbf{a}_0, \ldots,\mathbf{a}_{j_{k}-1},b_k,\mathbf{a}_{j_k},\ldots, \mathbf{a}_{j})= \tilde G \cdot (\mathbf{a}_{j_k},\ldots, \mathbf{a}_{j})$, for some $k \geq 1$ and $j_k\leq j<j_{k+1}$, then define
\[\lambda(G) = \lambda(\tilde G^*)  \prod_{l=0}^{j-j_k} \bar\nu_p(\mathbf{a}_{j_{k}+l}).\]
\end{itemize}
Recall that we refer to those sequences obtained as above as admissible sequences. These describe the support of our measure. With this definition in hand, we verify the continuant growth estimate \eqref{typicalcontinuant}.
\begin{proof}[Proof of inequality \eqref{typicalcontinuant}]
The proof of this will follow from inductively applying the gluing lemma, Lemma \ref{gluinglemma}. We will show the following statements using induction. 
\begin{itemize}
\item $\heartsuit(k)$ for $k \geq 0$: If $G$ is of the form $G = (\mathbf{a}_0, \ldots, \mathbf{a}_{j_k - 1} b_k, \mathbf{a}_{j_k}, \ldots, \mathbf{a}_{j-1})$ where $j_k \leq j < j_{k+1}$ (with the understanding that $j_0 = 0$), then
\[\left| \log K(G) - (j + (\tau - 2) j_k) \sigma \right| < \epsilon j \sigma.\]
If $j = j_{k+1}$, then the bound $\epsilon j \sigma$ can be improved to $\frac{\epsilon}{100} j \sigma$.
\item $\spadesuit(k)$ for $k \geq 1$: If $G$ is of the form $G = (\mathbf{a}_0, \ldots \mathbf{a}_{j_k - 1}, b_k)$, then 
\[\left| \log K(G) - (j_k + (\tau - 2) j_k) \sigma \right| < \frac{1}{2} \epsilon j_k \sigma.\]
\end{itemize}
\paragraph*{Proof of $\heartsuit(0)$} Let $G = (\mathbf{a}_0, \ldots, \mathbf{a}_{j-1})$ be such that $j < j_1$. By applying Lemma \ref{gluinglemma} $j - 1$ times, we arrive at the estimate
\[\left| \log K(G) -  \sum_{i=0}^{j-1} \log K(\mathbf{a}_i) \right| \leq (j-1) C_N\]
and
\[\left| \sum_{i=0}^{j-1} \log K (\mathbf{a}_i) - j \sigma \right| \leq \frac{\epsilon}{500} j \sigma\]

If $p$ is sufficiently large to ensure $\sigma \geq 500\epsilon^{-1}C_N$, then the triangle inequality gives
\begin{equation}\label{eq:KGgrowthOnlyGood}\left| \log K(G) - j \sigma \right| \leq \frac{\epsilon}{100} j \sigma\end{equation}
as desired.
\paragraph*{Proof of $\spadesuit(k)$ given $\heartsuit(k-1)$}
Suppose $G = (\mathbf{a}_0, \ldots, \mathbf{a}_{j_k - 1}, b_k)$. Let $G' = (\mathbf{a}_0, \ldots, \mathbf{a}_{j_k - 1})$. The inductive assumption $\heartsuit(k-1)$ implies that
\[\left|\log K(G') - (j_k + (\tau - 2) j_{k-1}) \sigma \right| < \frac{\epsilon}{100} j_k \sigma.\] 
If $j_k$ is chosen sufficiently large that $(\tau - 2) j_{k - 1} \sigma < \frac{\epsilon}{100} j_k \sigma$, then we have
\[\left|\log K(G') - j_k \sigma \right| < \frac{\epsilon}{50} j_k \sigma.\] 
By \eqref{rhoorder}, we can select a $j_k$ sufficiently large such that
\[\left|\frac{\log \rho(q)}{\log q} - (\tau - 2) \right| < \frac{\epsilon}{100},\]
provided that $q > e^{j_k \sigma / 2}$. Of course, $K(G')$ is significantly larger than $e^{j_k \sigma/2}$, so we have from \eqref{eq:bkInt} that 
\begin{IEEEeqnarray*}{Cl}
& \left| \log b_k - (\tau - 2) j_k \sigma \right| \\
\leq & \left| \log b_k - \log \rho(K(G')) \right|  + |\log \rho(K(G')) - (\tau - 2)\log K (G')| + \left|(\tau - 2)\log K(G') - (\tau - 2) j_k \sigma \right|\\
\leq & \log \left(1 + \frac{\eta_k}{50} \right) + \frac{\epsilon}{100} \log K(G) + (\tau - 2) \frac{\epsilon}{50} j_k \sigma \\
\leq & \log \left(1 + \frac{\eta_k}{50} \right) + \frac{\epsilon}{50} j_k \sigma + \frac{\epsilon}{50} j_k \sigma \\
\leq & \frac{\epsilon}{20} j_k \sigma,
\end{IEEEeqnarray*}
provided that $j_k$ is sufficiently large. Finally, we observe that $K(G) = b_k K(G') + K'(G')$, which implies $| \log K(G) - (\log b_k + \log K(G')) | \leq 1$. Therefore,
\begin{IEEEeqnarray*}{Cl}
& \left|\log K(G) - (j_k +  (\tau - 2) j_k) \sigma \right| \\
\leq & \left| \log K(G) - (\log b_k + \log K(G')) \right| + \left| \log b_k - (\tau - 2) j_k \sigma \right| + \left| \log K(G') - j_k \sigma \right| \\
\leq & 1 + \frac{\epsilon}{20} j_k \sigma + \frac{\epsilon}{50} j_k \sigma \\
\leq & \frac{\epsilon}{2} j_k \sigma,
\end{IEEEeqnarray*}
as desired.
\paragraph*{Proof of $\heartsuit(k)$ given $\spadesuit(k)$ for $k \geq 1$} Let $G = (\mathbf{a}_0, \ldots, \mathbf{a}_{j_k - 1}, b_k, \mathbf{a}_{j_k}, \ldots, \mathbf{a}_{j-1})$ where $j_k < j \leq j_{k+1}$. Let $G' = (\mathbf{a}_0, \ldots, \mathbf{a}_{j_k - 1}, b_k)$. By applying Lemma \ref{gluinglemma} a total of $j - j_k$ times, we get
\[\left| \log K(G) - \log K(G') - \sum_{i = j_k}^{j - 1} \log K(\mathbf{a}_i) \right| \leq C_N (j - j_k)\]
But
\[\left| \log K(G') - (j_k + (\tau - 2) j_k) \sigma \right| < \frac{\epsilon}{2} j_k \sigma\]
and
\[\left|\sum_{i=j_k}^{j-1} \log K(a_i) - (j - j_k) \sigma \right| < \frac{\epsilon}{500} (j - j_k) \sigma < \frac{\epsilon}{500} j \sigma.\]
Hence
\begin{IEEEeqnarray*}{CCl}
& & \left| \log K(G) - (j + (\tau - 2) j_k) \sigma \right| \\
\leq & & \left| \log K(G) - \log K(G') - \sum_{i=j_k}^{j-1} \log K(\mathbf{a}_i)\right| \\
& + & \left| \log K(G') - (j_k + (\tau - 2) j_k) \sigma \right| \\
& + & \left| \sum_{i=j_k}^{j - 1}\log K(\mathbf{a}_i) - (j - j_k) \sigma \right|  \\
\leq & & C_N (j - j_k) + \frac{\epsilon}{2} j_k \sigma + \frac{\epsilon}{500} j \sigma.
\end{IEEEeqnarray*}
It is then easy to see that this is bounded above by $\epsilon j \sigma$. Indeed, this follows because $j > j_k$ and $\sigma \gg C_N$. If $j = j_{k+1}$ and $j_{k+1}$ is sufficiently large depending on $j_k$, then we have the improved bound  $\frac{\epsilon}{100} j \sigma$.\end{proof}

\section{Geometry of the measure}\label{sec:measureOfBalls}
\setcounter{subsection}{1}
In this section, we derive some of the key geometric properties of our measure. Specifically, these are upper and lower bounds on the $\lambda^\sharp$ measure of balls.

The following lemma is a statement about mass distributions arising from our measure construction. We call this the capping lemma because in practice, we will cap off a long sequence of small partial quotients with a large partial quotient.

\begin{mylem}[Capping lemma]\label{cappinglemma}
Let $\mu$ be a measure and let $G := (c_0, c_1, \ldots, c_n)$, be an admissible sequence, with $n$ sufficiently large. Suppose that $K(c_0, \ldots, c_n) = q_n$ and $K(c_0, \ldots, c_{n-1}) = q_{n-1}$ and suppose that $\mu(\cyl(G)) \leq  q_n^{-2 + O(\epsilon)} = |\cyl(G)|^{1 - O( \epsilon)}$. Suppose that the mass $\mu(\cyl(G))$ is distributed evenly among those subsets $\cyl(G \cdot b)$ for some $b \in [ \left(1 + \frac{\eta}{100} \right) \rho(q_n), \left( 1 + \frac{\eta}{50} \right) \rho(q_n)]$. Then 
\begin{equation}\label{muJestimate}
\mu(\cyl(G \cdot b)) \leq \eta^{-1} |\cyl(G \cdot b)|^{\frac{\tau }{2 \tau - 2} - O(\epsilon)}.\end{equation} 
\end{mylem}
\begin{proof}
Suppose $\mu(\cyl(G)) \leq q^{-2 +O(\epsilon)}$. Let $\omega_n := \frac{\log \rho(q_n)}{\log q_n}\rightarrow \tau-2$ so that, for $q_n$ sufficiently large (which we can ensure by taking $n$ sufficiently large), we have the estimate $0<\omega_n < \tau - 2 + \epsilon$. 
Then the mass $\mu$ is distributed evenly across $\frac{\eta}{1000} \rho(q_n) = \frac{\eta}{1000} q_n^{\omega_n}$ subsets $\cyl(G \cdot b)$ of $\cyl(G)$. So each such subset $\cyl(G \cdot b)$ will have $\mu$-measure at most
\[\mu(\cyl(G \cdot b)) \leq 1000 \eta^{-1} q^{-2  - \omega_n+O(\epsilon)}.\]
However, the length of $\cyl(G \cdot b)$ satisfies $|\cyl(G \cdot b)| \sim K(c_0, \ldots, c_m, b)^{-2}$ and $b \sim \rho(q_n) = q_n^{\omega_n}$, so 
\begin{IEEEeqnarray*}{rCl}
K(c_0, \ldots, c_n, b)^{-2} & = & (bq_n + q_{n-1})^{-2} \\
& \geq & ((b + 1) q_n)^{-2} \\
& \geq & \frac{1}{2} (q_n^{1 + \omega_n})^{-2} \\
& = & \frac{1}{2} q_n^{-2 - 2 \omega_n}.
\end{IEEEeqnarray*}
So, for an appropriate constant $C$, we have
\begin{IEEEeqnarray*}{rCl}
\mu(\cyl(G \cdot b)) & \leq & C_N' \eta^{-1} |\cyl(G \cdot b)|^{\frac{ 2  + \omega_n - O(\epsilon)}{2 + 2 \omega_n}} \\
& \leq &  \eta^{-1} |\cyl(G \cdot b)|^{\frac{1}{2}+\frac{1 }{2 + 2 (\tau - 2 )} - O(\epsilon)}.
\end{IEEEeqnarray*}
\end{proof}

In order to discuss the measure constructed in Section \ref{Measure}, it is useful to be able to refer to particular finite sequences. Recall that a sequence is called admissible if it is made up of $\mathbf{a}_j$ lying in the support of $\bar \nu_p$ and suitably separated entries $b_k$ satisfying \eqref{eq:bkInt}.
\begin{mydef}
An admissible sequence of the form 
\[(\mathbf{a}_0, \ldots, \mathbf{a}_{j_1 - 1}, b_1, \mathbf{a}_{j_1}, \ldots, \mathbf{a}_{j_2 - 1}, b_2, \ldots, \mathbf{a}_j)\]
 will be called an $a$-sequence. An admissible sequence of the form 
\[(\mathbf{a}_0, \ldots, \mathbf{a}_{j_1 - 1}, b_1, \mathbf{a}_{j_1}, \ldots, \mathbf{a}_{j_2 - 1}, b_2, \ldots, \mathbf{a}_{j_k - 1}, b_k)\]
is called a $b$-sequence.

An $a$-sequence will be called \textbf{terminal} if it terminates in $\mathbf{a}_{j_k - 1}$ for some $j_k$; otherwise it will be called \textbf{non-terminal}. Thus a $b$-sequence or a non-terminal $a$-sequence will always be followed by a $p$-tuple of typical partial quotients $\mathbf{a}$, and a terminal $a$-sequence will always be followed by an exceptional partial quotient $b$.
\end{mydef}
We will show the following statement about the measure $\lambda^\sharp$.
\begin{mylem}\label{cylinderlemma}
Let $\{\eta_k\}_{k=1}^{\infty}$ be a sequence that decays to zero (for example, $\eta_k = 2^{-k}$ would work). Then provided that the sequence $\{j_k\}_{k=1}^{\infty}$ grows sufficiently quickly depending on $\epsilon$, we have the following estimates on $\lambda$:
\begin{itemize}
\item If $G = (\mathbf{a}_0, \ldots, \mathbf{a}_{j_k - 1})$ is a terminal $a$-sequence, then $\lambda^\sharp(\cyl(G)) \leq |\cyl(G)|^{1 - O( \epsilon)}$.
\item If $G = (\mathbf{a}_0, \ldots, \mathbf{a}_{j_k - 1}, b_k)$ is a terminal $b$-sequence, then $\lambda^\sharp(\cyl(G)) \leq |\cyl(G)|^{\frac{\tau }{2 \tau - 2} - O(\epsilon)}$
\item If $G = (\mathbf{a}_0, \ldots, \mathbf{a}_{j})$ is an $a$-sequence, then $\lambda^\sharp(\cyl(G)) \leq  |\cyl(G)|^{\frac{\tau }{2 \tau - 2} - O(\epsilon)}$.
\end{itemize}
\end{mylem}
\begin{proof}
We first establish the result for terminal sequences by induction and then show how the remaining bound for general $a$-sequences follows.
\paragraph*{Case 1: $G$ is the first terminal $a$-sequence.} We will start with the case $G = (\mathbf{a}_0, \ldots, \mathbf{a}_{j_1 - 1})$. 

We observe that
\begin{IEEEeqnarray*}{Cl}
& \bar \nu_p(a_0, \ldots, a_{p-1}) \cdots \bar \nu_p(a_{(j_1-1)p}, \ldots, a_{j_1 p - 1}) \\
\leq & 2^{j_1} \nu_m(a_0, \ldots, a_{m-1}) \cdots \nu_m(a_{(j_1-1)p + (J - 1)m} ,\ldots, a_{j_1p - 1}) \\
= & 2^{j_1} S_m^{-j_1 J} K(a_0, \ldots, a_{m-1})^{-2 (1 - \epsilon)} \cdots K(a_{(j_1-1) p + (J - 1)m}, \ldots, a_{j_1 p - 1})^{-2 (1 - \epsilon)} \\
\leq & 2^{j_1} S_m^{-j_1 J} e^{C_N (j_1 J - 1)} K(a_0, \ldots, a_{j_1 p - 1})^{-2 (1 - \epsilon)} \\
\leq & 2^{j_1} C_N'S_m^{-j_1 J} e^{C_N (j_1 J - 1)} |\cyl(G)|^{1 - \epsilon}
\end{IEEEeqnarray*}
Here, the last inequality is obtained by applying Lemma \ref{gluinglemma}, the gluing lemma, a total of $j_1 J - 1$ times. If $m$ is chosen sufficiently large that $S_m > C_N'2e^{C_N}$, then this expression will be bounded above by $|\cyl(G)|^{1 - \epsilon}$.

\paragraph*{Case 2: $G$ is a $b$-sequence.} Second, let $G = (\mathbf{a}_0, \ldots, \mathbf{a}_{j_k - 1}, b_k)$. Now, assume by induction that we know that $\lambda^\sharp(\cyl(\tilde G)) \leq |\cyl(\tilde G)|^{1 - \epsilon}$, where $\tilde G = (\mathbf{a}_0, \ldots, \mathbf{a}_{j_k - 1})$. We will use Lemma \ref{cappinglemma}, the capping lemma, to show that $\lambda^\sharp(\cyl(G)) \leq |\cyl(G)|^{\frac{\tau }{2 \tau - 2}-O(\epsilon)}$ for $G = (\mathbf{a}_0, \ldots, \mathbf{a}_{j_k - 1}, b_k)$ provided that $j_k$ is sufficiently large. In fact, the capping lemma gives the bound
\[\lambda^\sharp(\cyl(G)) \leq  \eta_k^{-1} |\cyl(G)|^{\frac{\tau }{2 \tau - 2} - O(\epsilon)}\]
So the desired bound will hold, provided that $j_k$ is chosen sufficiently large to ensure $|\cyl(G)|^{\frac{\epsilon}{2\tau - 2}} < \eta_k$.

\paragraph*{Case 3: $G$ is another terminal $a$-sequence.} Finally, assume $G = (\mathbf{a}_0, \ldots, \mathbf{a}_{j_k - 1})$ for some $k \geq 2$. Then take $\tilde G$ to be $(a_0, \ldots, a_{j_{k-1} p - 1}, b_{k-1})$. 
By definition of $\lambda^\sharp(\cyl(G))$, we have
\[\lambda^\sharp(\cyl(G)) = \lambda^\sharp(\cyl(\tilde G)) \bar \nu_p(a_{j_{k-1} p}, \ldots, a_{j_{k-1}p + p - 1}) \cdots \bar \nu_p(a_{(j_k - 1) p}, \ldots a_{j_k p - 1}).\]
We know from Case 2 that
\[\lambda^\sharp(\cyl(\tilde G))\leq   |\cyl(\tilde G)|^{\frac{\tau }{2 \tau - 2} - O(\epsilon)}=K(\tilde G)^{\frac{-2\tau }{2 \tau - 2} + O(\epsilon)} = K(\tilde G)^{-2(1-\epsilon)+2(1-\epsilon) -\frac{2\tau }{2 \tau - 2}+O(\epsilon)}.\]
Provided we choose $j_k$ sufficiently large, this factor will be insubstantial.
The rest of the estimate proceeds as in Case 1:
\begin{IEEEeqnarray*}{Cl}
& \bar \nu_p(a_{j_{k-1} p}, \ldots, a_{j_{k-1}p + p - 1}) \cdots \bar \nu_p(a_{(j_k - 1) p}, \ldots, a_{j_k p - 1}) \\
\leq & 2^{j_k - j_{k-1}} S_m^{-(j_k - j_{k-1}) J } K(a_{j_{k-1} p }, \ldots, a_{j_{k-1} p + m-1})^{- 2 (1 - \epsilon)} \cdots K(a_{j_k - p + (J - 1)m}, \ldots, a_{j_k - 1})^{-2 (1 - \epsilon)} \\
\leq & 2^{j_k - j_{k-1}} S_m^{-(j_k - j_{k-1}) J} e^{C_N ((j_k - j_{k-1}) J - 1)} K(a_{j_{k-1}p}, \ldots, a_{j_k p - 1})^{-2 (1 - \epsilon)},
\end{IEEEeqnarray*}
where the last step follows from applying Lemma \ref{gluinglemma} $(j_k - j_{k-1}) J - 1$ times. Applying the gluing lemma once more, we find that
\begin{IEEEeqnarray*}{Cl}
\lambda^\sharp(\cyl(G)) &\leq  \left(K(\tilde G)^{2(1-\epsilon) -\frac{2\tau }{2 \tau - 2}+O(\epsilon)}\right)\left(2^{j_k - j_{k-1}} S_m^{-(j_k - j_{k-1}) J} e^{C_N ((j_k - j_{k-1}) J)} K(G)^{-2 (1 - \epsilon)}\right) \\
& \leq \left(C_N'K(\tilde G)^{2(1-\epsilon) -\frac{2\tau }{2 \tau - 2}+O(\epsilon)}2^{j_k - j_{k-1}} S_m^{-(j_k - j_{k-1}) J} e^{C_N ((j_k - j_{k-1}) J)} \right)\left(|\cyl(G)|^{1 - \epsilon}\right).
\end{IEEEeqnarray*}
If $m$ is chosen sufficiently large that $S_m > 2e^{C_N}$, then the left factor will be bounded above by $1$, provided that $j_k - j_{k-1}$ is sufficiently large. 

We have established the result for terminal sequences. For a general $a$-sequence, we can make use of what is known about the $b$-sequences. Assume that $G = (\mathbf{a}_0, \ldots, b_k, \mathbf{a}_{j_k },\ldots, \mathbf{a}_{j})$, with $j\leq j_{k+1}-1$. We write $G=\tilde{G}\cdot H$, where $\tilde{G}=(\mathbf{a}_0, \ldots, b_k)$ and $H=(\mathbf{a}_{j_k },\ldots, \mathbf{a}_{j})$. By definition, what we have already shown, and the gluing lemma,
\begin{IEEEeqnarray*}{Cl}\lambda^\sharp(\cyl(G)) &= \lambda^\sharp(\cyl(\tilde{G}))\bar\nu_p(\mathbf{a}_{j_k })\ldots \bar\nu_p(\mathbf{a}_{j})\\
& \leq  |\cyl(\tilde{G})|^{\frac{\tau }{2 \tau - 2} - O(\epsilon)}|\cyl(H)|^{1-\epsilon}\\
& \leq  |\cyl(\tilde{G})|^{\frac{\tau }{2 \tau - 2} - O(\epsilon)}|\cyl(H)|^{\frac{\tau }{2 \tau - 2} -\epsilon}\\
& \leq  |\cyl(G)|^{\frac{\tau }{2 \tau - 2} - O(\epsilon)}.
\end{IEEEeqnarray*}
\end{proof}
We have established the behavior of the measure $\lambda^\sharp$ on cylinder sets, from which the ball condition for $\lambda^{\sharp}$ will follow.
\begin{mypro}\label{prop:unifBallCond}
Let $I = [t, t + h]$ be any interval. Then, for sufficiently small $h$, $\lambda^{\sharp}(I) \leq h^{\frac{2}{\tau} - O(\epsilon)}$ provided that the $j_k$ grow quickly enough.
\end{mypro}
\begin{proof}
If the interval $I$ intersects the support of $\lambda^{\sharp}$, then there is a longest $a$ or $b$ sequence $G$ such that $I \cap \supp \lambda^{\sharp} \subset \cyl(G)$. We will split into two cases based on whether the sequence $G$ is a terminal $a$-sequence. 
\paragraph*{Case 1: $G$ is either a non-terminal $a$-sequence or a $b$-sequence} If $G$ is not a terminal $a$-sequence, then $G$ is either of of the form $(\mathbf{a}_0, \ldots, \mathbf{a}_{r})$ for some $r \neq j_k - 1$ or of the form $(\mathbf{a}_0, \ldots, \mathbf{a}_{j_k - 1}, b_k)$. We claim that $|\cyl(G)|\leq A_N h$ for some constant $A_N$ depending only on $N$. Indeed, none of the cylinder sets $\cyl(G')$ where $G'$ is an $a$-sequence of the form $(\mathbf{a}_0, \ldots, \mathbf{a}_{r + 1})$ or $(\mathbf{a}_0, \ldots, \mathbf{a}_{j_k - 1}, b_k, \mathbf{a}_{j_k})$ contain all of $[t, t + h] \cap \supp \lambda^{\sharp}$. This means that $h$ must be at least the minimum of the distances between such intervals $\cyl(G')$. But the distance between such intervals is at least an $N,p$-dependent constant times $K(G)^{-2} \gtrsim_{N,p} |\cyl(G)|$. Thus the desired upper bound on $\lambda^{\sharp}(I)$ follows from the cylinder set estimate from Lemma \ref{cylinderlemma}.
\paragraph*{Case 2: $G$ is a terminal $a$-sequence} We now consider the complementary case in which $G$ is of the form $(\mathbf{a}_0, \ldots, \mathbf{a}_{j_k - 1})$. For such values of $G$, there is no convenient estimate available on $h$. We split into two subcases based on the size of $h$. We will write 
\[\omega := \frac{\log \rho(K(G))}{\log K(G)}\]
so that $\rho(K(G)) = K(G)^{\omega}$.  For the rest of this argument, we will assume $h$ is sufficiently small to ensure, via \eqref{rhoorder}, that $\omega \leq \tau - 2 + \epsilon$.
\subparagraph*{Case 2a} We first consider the case in which $h > |\cyl(G)|^{(\omega + 2)/2 }$. By Lemma \ref{cylinderlemma}, we have that 
\begin{IEEEeqnarray*}{Cl}
& \lambda^{\sharp}(I)  \\
\leq & \lambda^\sharp(\cyl(G)) \\ 
\leq & |\cyl(G)|^{1 - O(\epsilon)} \\
\leq & h^{(1 - O(\epsilon)) \left(\frac{1}{(\omega + 2)/2} \right)} \\
\leq & h^{2/\tau - O(\epsilon)}.
\end{IEEEeqnarray*}
\subparagraph*{Case 2b} We next consider the case in which $h \leq |\cyl(G)|^{(\omega + 2)/2}$. We use $G\cdot b$ to denote an admissible tuple $(\mathbf{a}_0, \ldots, \mathbf{a}_{j_k - 1}, b)$. We know that $b\in [(1+\eta_k/100)K(G)^{\omega},(1+\eta_k/50)K(G)^{\omega}]$. The number of choices of $b$ is at least $C_N' \eta_k |\cyl(G)|^{-\omega/2}$ (here, $C_N'$ denotes a constant which may change from line to line). Since the measure of $\cyl(G)$ is distributed evenly among the $\cyl({G \cdot b)}$, using the inductive assumption we have that 
\begin{IEEEeqnarray*}{rCl}
\lambda^\sharp(\cyl({G \cdot b})) & \leq & C_N' \eta_k^{-1} |\cyl(G)|^{\omega/2} \lambda^\sharp(\cyl(G))  \\
& \leq &  \eta_k^{-1} |\cyl(G)|^{1 + \omega/2 - O(\epsilon)}.
\end{IEEEeqnarray*}
Furthermore, the length of the interval $\cyl(G \cdot b)$ is at least $C_N' |\cyl(G)|^{(1 + \omega)}$. So the pigeonhole principle implies that the number of such intervals $\cyl(G \cdot b)$ intersecting $I$ is no more than $h C_N' |\cyl(G)|^{- \omega - 1}$. Using the cylinder estimates, we thus find that 
\begin{IEEEeqnarray*}{rCl}
\lambda^\sharp(I) & \leq &  h |\cyl(G)|^{- \omega - 1} \eta_k^{-1} |\cyl(G)|^{1 + \omega/2 - O(\epsilon)} \\
& \leq &  \eta_k^{-1} h |\cyl(G)|^{- \omega/2 -  O(\epsilon)} \\
& \leq & \eta_k^{-1} h \cdot h^{\frac{2}{\omega + 2} (- \omega/2 -  O(\epsilon))} \\
& \leq & \eta_k^{-1} h^{\frac{2}{\omega + 2} - O(\epsilon)}.
\end{IEEEeqnarray*}
But $\omega + 2 \leq \tau + \epsilon$, so provided that $j_k$ is chosen sufficiently large depending on $\eta_k$, we have that the above expression is at most $h^{\frac{2}{\tau} - O(\epsilon)}$.
\end{proof}

As we are going to be making use of \textit{normalized} pushforward measures, it is also necessary to record the following \textit{lower} bound. 

\begin{mypro}\label{cylinderlowerbound}
Let $G$ be an admissible sequence. Then 
\begin{equation}\label{lambdalowerbound}
\lambda^\sharp(\cyl(G)) \gtrsim |\cyl(G)|^{1-O(\epsilon)}.
\end{equation}
\end{mypro}
Before proceeding with the proof of Proposition \ref{cylinderlowerbound}, we require some preliminary results. Recall that for a large integer $m$, we defined $\nu_m$ to be a measure on $\{1, 2, \ldots, N\}^m$ defined by
\begin{equation}\label{numdef}
\nu_m(a_0, \ldots, a_{m-1}) = \frac{1}{S_m} K(a_0, \ldots, a_{m-1})^{-2 (1-\epsilon)}.
\end{equation}
For this argument the precise value of $S_m$ will not be of any consequence. Recall also that $p = Jm$ for some large integer $J$,  and $\bar \nu_p$ is a measure on $\{1, 2, \ldots, N\}^m$ satisfying the condition
\begin{equation}\label{nupcon}
\prod_{i=1}^{J} \nu_m(a_{(i-1)m}, \ldots, a_{im - 1}) \leq \bar \nu_p(a_0, \ldots, a_{p-1}) \leq 2 \prod_{i=1}^{J} \nu_m(a_{(i-1)m}, \ldots, a_{im - 1}).
\end{equation}

Recall that we chose $N$ to ensure that
\[1 - \epsilon < \text{dim}_H (\textbf{Bad}(N)) < 1 - \epsilon/10.\]
We are now ready to state the following lemma.
\begin{mylem}
Let $G_1, G_2 \in \{1, 2, \cdots, N\}^{jp}$ be $jp$-tuples in the support of $(\bar \nu_p)^j$, with $j$ sufficiently large. Then 
\begin{equation}\label{lebesguecomparison}
|\cyl (G_1)|^{1 + \epsilon/50}  \leq |\cyl (G_2)| \leq |\cyl (G_1)|^{1 - \epsilon/50} 
\end{equation}
and
\begin{equation}\label{nupcomparison}
(\bar \nu_p \times \cdots \times \bar \nu_p(G_1))^{1+\epsilon/50}\leq \bar \nu_p \times \cdots \times \bar \nu_p(G_2)\leq(\bar \nu_p \times \cdots \times \bar \nu_p(G_1))^{1 - \epsilon/50}.
\end{equation}
\end{mylem}
This lemma states that for all $jp$-tuple in the support of $(\bar \nu_p)^j$, the Lebesgue measures are comparable, and the $(\nu_p)^j$ measures are comparable. This will mean in practice that the \textit{ball conditions} satisfied by $\cyl(G_1)$ and $\cyl(G_2)$ will be comparable.
\begin{proof}
Let $G_1 = (a_0, \ldots, a_{jp - 1})$ and $G_2 = (a_0', \ldots, a_{jp - 1}')$.
For $G_1, G_2 \in \supp (\nu_p)^j$, we have the estimate
\begin{equation}\label{eq:KG1G2comp}(1 - \epsilon/100) j \sigma \leq \log K(G_1), \log K(G_2) \leq (1 + \epsilon/100) j \sigma,\end{equation}
like we showed in \eqref{eq:KGgrowthOnlyGood}. Since $|\cyl(G_1)| \sim |K(G_1)|^{-2}$ and $|\cyl(G_2)| \sim |K(G_2)|^{-2}$, we get the estimate \eqref{lebesguecomparison}.

It remains to prove \eqref{nupcomparison}. Because of the product structure of the measure, is suffices to prove it for sequences $G_1$ and $G_2$ in the support of $\bar\nu_p$. 
Using \eqref{nupcon} and \eqref{numdef}, we have the following estimate:
\begin{IEEEeqnarray*}{Cl}
& \bar \nu_p (a_0, \ldots, a_{p - 1}) \\
\geq  & \nu_m(a_0, \ldots, a_{m-1}) \times \cdots \times \nu_m(a_{p - m}, \ldots, a_{p - 1}) \\
= & \left(\frac{1}{S_m} \right)^{J} K(a_0, \cdots a_{m-1})^{-2 (1-\epsilon)} \cdots K(a_{p - m}, \ldots, a_{p - 1})^{-2 (1-\epsilon)} \\
\geq & e^{-C_N (J  - 1)} \left(\frac{1}{S_m} \right)^{J } K(a_0, \ldots, a_{p - 1})^{-2 (1-\epsilon)}.
\end{IEEEeqnarray*}
Similarly, we have the upper bound
\begin{IEEEeqnarray*}{Cl}
& \bar \nu_p (a_0, \ldots, a_{p - 1}) \\
\leq & 2 e^{C_N (J  - 1)} \left(\frac{1}{S_m} \right)^{J } K(a_0, \ldots, a_{p - 1})^{-2 (1-\epsilon)}.
\end{IEEEeqnarray*}
In this way, we find that 
\begin{IEEEeqnarray*}{Cl}
& \frac{\bar \nu_p(G_1)^{1 + \epsilon/50}}{\bar \nu_p(G_2)} \\
\leq &\left( \frac{2 e^{C_N (J  - 1)} K(G_{1})^{-2(1-\epsilon)}}{S_m^J}\right)^{1+\epsilon/50}\left( \frac{e^{C_N (J  - 1)}S_m^J}{K(G_{2})^{-2(1-\epsilon)}}\right).\end{IEEEeqnarray*}
We know by construction that $K(G_1),K(G_2)\in [e^{\sigma-\sigma\epsilon/500},e^{\sigma+\sigma\epsilon/500}]$ so that $K(G_1)\geq K(G_2)^{\frac{1-\epsilon/500}{1+\epsilon/500}}\geq K(G_2)^{1-\epsilon/200}$. Thus we see that
\begin{IEEEeqnarray*}{Cl}
& \frac{\bar \nu_p(G_1)^{1 + \epsilon/50}}{\bar \nu_p(G_2)} \\
\leq & 2^{1+\epsilon/50} e^{C_N (J  - 1)(2+\epsilon/50)} K(G_{2})^{-2(1-\epsilon)\left((1-\epsilon/200)(1+\epsilon/50)-1\right)}\left( \frac{1}{S_m^{J\epsilon/50}}\right)\\
\leq & 1,
\end{IEEEeqnarray*}
provided our initial choice of $m$ is sufficiently large. 
Similarly, we can show that 
\begin{IEEEeqnarray*}{Cl}
& \frac{\bar \nu_p(G_2)} {\bar\nu_p(G_1)^{1 - \epsilon/50} } \\
\leq & e^{C_N (J  - 1)(2+\epsilon/50)} K(G_{1})^{-2(1-\epsilon)\left((1-\epsilon/200)-(1-\epsilon/50)\right)}\left( \frac{1}{S_m^{J\epsilon/50}}\right)\\
\leq & 1.
\end{IEEEeqnarray*}
\end{proof}

The comparability of $(\bar \nu_p)^j(G_1)$ and $(\bar \nu_p)^j(G_2)$ is enough to give a lower bound on $(\bar \nu_p)^j(G)$ for all $G$.
\begin{mylem}\label{nupjballconlemma}
For any $G \in \supp (\bar \nu_p)^j$, we have 
\begin{equation}\label{nupjballcon}
(\bar \nu_p)^j(G) \gtrsim |\cyl(G)|^{1 -\epsilon/1000}.
\end{equation}
In particular, the implicit constant does not depend on $j$.
\end{mylem}
\begin{proof}
Suppose the estimate does not hold. Then for every $j$, there exists a set $G_0(j) \in \supp(\bar \nu_p)^j$ such that
\begin{equation}\label{nupjballconwrongassn}
(\bar \nu_p)^j(G_0) \lesssim |\cyl(G_0)|^{1 - \epsilon/1000}.
\end{equation}
In fact, equation \eqref{nupjballconwrongassn} implies an estimate on $(\bar \nu_p)^j(G)$ for \textit{all} $G \in \supp (\bar \nu_p)^j$. Also using the comparisons \eqref{lebesguecomparison} and  \eqref{nupcomparison}, we have 
\[(\bar \nu_p )^j(G) \lesssim |G|^{\frac{(1-\epsilon/1000)(1-\epsilon/50)}{(1+\epsilon/50)}}\]
so for all $G \in \supp(\bar \nu_p)^j$, we have
\begin{equation}\label{nupjballconwrongassn2}
(\bar \nu_p )^j(G) \lesssim |\cyl(G)|^{(1 - \epsilon/15)}
\end{equation}
Now, consider the measure $\bar \nu_p \times \bar \nu_p \times \cdots$ on $\{1, \ldots, N\}^{\infty}$. This measure lifts to a Borel measure $\nu^{\sharp}$ on $\mathbb{R}$  under the continued fraction map. By \eqref{nupjballconwrongassn2}, we have for any $j$ and any set $G(j)$ in the support of $(\bar \nu_p)^j$ that 
\[(\bar \nu_p)^j(G) \lesssim |\cyl(G)|^{1 - \epsilon/15}.\] 
Now, let $h > 0$ be a small number and consider an arbitrary interval $[t, t + h] \subset \mathbb{R}$. We will use the term \textbf{$p$-cylinder set} to refer to a set of the form $\cyl(G)$, where $G \in \{1, \ldots, N\}^{jp}$ for some integer $j \geq 0$. Let $G = (a_0, \ldots, a_{jp - 1})$ be such that $\cyl(G)$ is the minimal $p$-cylinder set containing $[t, t +h] \cap \text{Bad}(N)$. Because $\cyl(G)$ is the minimal such $p$-cylinder set, it follows that if $G' = (a_0, \ldots, a_{(j+1)p - 1})$, then $\cyl(G') \cap \text{Bad}(N)$ cannot contain all of $[t, t + h] \cap \text{Bad}(N)$.  Therefore, the length $h$ must be at least as large as the gaps between such sets $\cyl(G')$; the length of each such gap is $\gtrsim_{p, N} K(G)^{-2} \gtrsim |\cyl(G)|$. Hence we have
\begin{IEEEeqnarray*}{Cl}
& \nu^{\sharp}([t, t + h]) \\
\leq & \nu^{\sharp}(\cyl(G)) \\
= & (\bar \nu_p)^j(G) \\
\lesssim & |\cyl(G)|^{1 - \epsilon/15} \\
\lesssim & h^{1 - \epsilon/15}.
\end{IEEEeqnarray*}
Hence $\nu^{\sharp}$ is a Borel probability measure on $\text{Bad}(N)$ satisfying a $1 - \epsilon/15$ ball condition. But $N$ was chosen so that $\text{Bad}(N)$ would have Hausdorff dimension at most $1 - \epsilon/10$, so the existence of the measure $\nu^{\sharp}$ contradicts Frostman's lemma. The contradiction proves the estimate \eqref{nupjballcon}.
\end{proof}
With this lemma in hand, we are now ready to establish a lower bound on the size of all cylinder sets with respect to the measure $\lambda$.
\begin{proof}[Proof of Proposition \ref{cylinderlowerbound}]
Let $G = (\mathbf{a}_0, \ldots, \mathbf{a}_{j_1 - 1}, b_1, \mathbf{a}_{j_1}, \ldots, c)$ be a finite sequence, where $c$ is either a $p$-tuple of typical partial quotients or a single exceptional partial quotient. It is enough to show \eqref{lambdalowerbound} for such cylinder sets, as if $G$ ends with $a_t$ for $t$ not equal to $jp$, then both $|\cyl(G)|$ and $\lambda(\cyl(G))$ will be comparable to $|\cyl(G')|$ and $|\lambda(G')|$, where $G'$ is an appropriate substring of $G$ that ends with either some $b_k$ or $\mathbf{a}_j$ where $jp < t < (j + 1)p$.
We will show the following estimates using induction.
\begin{itemize}
\item $\clubsuit(k)$ for $k \geq 0$: If $G$ is of the form $G = (\mathbf{a}_0, \ldots, \mathbf{a}_{j_k - 1}, b_k, \mathbf{a}_{j_k}, \ldots, \mathbf{a}_j)$, where $j_k \leq j < j_{k+1}$ (or $0 \leq j < j_1$ if $k = 0$), then
\[\lambda^\sharp(\cyl(G)) \geq C_{N,p, \epsilon}^{2k+1} |\cyl(G)|^{1 - \epsilon/1000}\]
\item $\diamondsuit(k)$ for $k \geq 1$:  If $G$ is of the form $G = (\mathbf{a}_0, \ldots, \mathbf{a}_{j_k - 1}, b_k)$, then \[\lambda^\sharp(\cyl(G)) \geq C_{N,p, \epsilon}^{2k} |\cyl(G)|^{1 - \epsilon/1000}.\]
\end{itemize}
The estimate \eqref{lambdalowerbound} follows from the above induction by using the rapid growth rate of the indices $j_k$ to obtain the uniform estimate $C_{N,p, \epsilon}^{k}|\cyl(G)|^{ - \epsilon/2000}\gtrsim 1$, for corresponding $G$.
\paragraph*{Proof of $\clubsuit(0)$}
Since $G$ does not contain any exceptional partial quotients, the statement $\clubsuit(0)$ is in fact equivalent to the estimate $(\bar \nu_p)^j(G) \gtrsim_{N,p, \epsilon} |\cyl(G)|^{1 - \epsilon/100}$. This is precisely the estimate shown in Lemma \ref{nupjballconlemma}.
\paragraph*{Proof of $\diamondsuit(k)$ assuming $\clubsuit(k-1)$}
Assume $\clubsuit(k-1)$ holds for some $k$. Recall that $G$ is assumed to be a $b$-sequence; say $G := (\mathbf{a}_0, \ldots, \mathbf{a}_{j_k - 1}, b_k)$. Define $\tilde G = (\mathbf{a}_0, \ldots, \mathbf{a}_{j_k - 1})$ so that $\tilde G$ is a terminal $a$-sequence. The estimate $\clubsuit(k-1)$ implies 
\[\lambda^\sharp(\cyl(\tilde G)) \geq C_{N,p, \epsilon}^{2k - 1} |\cyl(\tilde G)|^{1 - \epsilon/1000}.\]
Write $\omega := \frac{\log \rho(K(G))}{\log K(G)}$. Since $K(G) \sim K(\tilde G) \cdot b_k \sim K(\tilde G) \cdot K(\tilde G)^{\omega}$, we have $K(G) \sim K(\tilde G)^{\omega + 1}$. Therefore, $|\cyl(G)| \sim |\cyl(\tilde G)|^{\omega + 1}$. On the other hand, we have by the definition of $\lambda^\sharp$ that 
\[\lambda^\sharp(\cyl(G)) \sim \lambda^\sharp(\cyl(\tilde G)) \cdot \eta_k^{-1} K(\tilde G)^{-\omega} \gtrsim \lambda^\sharp(\cyl(\tilde G)) \cdot |\cyl(\tilde G)|^{\omega/2}.\]
Therefore, if $\epsilon$ is small enough (depending on $\tau$) and $C_{N, p, \epsilon}$ is small enough, we have the following estimate from $\clubsuit(k-1)$
\begin{IEEEeqnarray*}{rCl}
\frac{\lambda^\sharp(\cyl(G))}{|\cyl(G)|^{1 - \epsilon/1000}} & \geq & C_{N, p, \epsilon}^{2k - 1} \frac{\lambda^\sharp(\cyl(\tilde G)) \cdot |\cyl(\tilde G)|^{\omega/2}}{|\cyl(\tilde G)|^{(\omega + 1)(1 - \epsilon/1000)}} \\
& = & C_{N, p, \epsilon}^{2k -1} \lambda^\sharp(\cyl(\tilde G)) |\cyl(\tilde G)|^{-\omega/2 - 1 + (\omega + 1) \epsilon/1000} \\
& \geq & C_{N, p, \epsilon}^{2k-1} \frac{\lambda^\sharp(\cyl (\tilde G))}{|\cyl(\tilde G)|^{1 - \epsilon/1000}},
\end{IEEEeqnarray*}
where the last step follows from the fact that $\omega \geq 0$, a consequence of \eqref{q2psiqlimit}.
\paragraph*{Proof of $\clubsuit(k)$ assuming $\diamondsuit(k)$, $k \geq 1$}
The remaining part of the induction is to establish $\clubsuit(k)$ from $\diamondsuit(k)$ for $k \geq 1$. By assumption, $G$ is of the form $(\mathbf{a}_0, \ldots, \mathbf{a}_{j_k - 1}, b_k, \mathbf{a}_{j_k}, \ldots, \mathbf{a}_j)$. We choose $\tilde G = (\mathbf{a}_0, \ldots, \mathbf{a}_{j_k - 1}, b_k)$ and $G' = (\mathbf{a}_{j_k}, \ldots, \mathbf{a}_{j})$. The definition of $\lambda$ implies that $\lambda(G) = \lambda(\tilde G^*) \cdot (\bar \nu_p)^{j - j_k + 1}(G')$. The inductive assumption $\diamondsuit(k)$ implies that $\lambda(\tilde G) \geq C_{N, p, \epsilon}^{2k} |\cyl(\tilde G)|^{1 - \epsilon/1000}$. Moreover, Lemma \ref{nupjballconlemma} implies that $(\bar \nu_p)^{j - j_k + 1}(G') \gtrsim |\cyl(G')|^{1 - \epsilon/1000}$. Hence, we have 
\[\lambda(G) \gtrsim C_{N, p, \epsilon}^{2k} |\cyl(\tilde G)|^{1 - \epsilon/1000} |\cyl(G')|^{1 - \epsilon/1000}.\]
Following the comparisons $K(\tilde G)^{-2}\sim|\cyl(\tilde G)|$ and $K( G')^{-2}\sim|\cyl( G')|$, we can apply the gluing lemma, Lemma \ref{gluinglemma}. Choosing $C_{N, p, \epsilon}^{2k+1}$ to be smaller than the implicit constant gives the desired bound.
\end{proof}

\section{Geometry of the relative measures}\label{sec:relMeasBalls}
\setcounter{subsection}{1}
Recall that, for two sequences $G$ and $F$, we write $G\cdot F$ for their concatenation. For any finite sequence $G$, and a collection of sequences $S$ 
, we defined the \textit{relative measure} $\lambda_G$ of $S$ by  
\[\lambda_{G}(S) := \frac{1}{\lambda(G^*)} \lambda(G\cdot S),\]
which defines a probability measure. We require ball conditions for the corresponding pushforward measures $\lambda_G^\sharp$.

We now make a quick observation that is important for understanding the behavior of the measures $\lambda_G$ and $\lambda_G^{\sharp}$. Let $g$ be the continued fraction map sending a finite or infinite sequence $F = (c_0, c_1, \ldots)$ to its continued fraction $[c_0; c_1, \ldots]$. The basic properties of continued fractions then guarantee that, for any $x \in \mathbb{R}$, we have
\begin{equation}\label{contfracshift}
g(G \cdot g^{-1}(x)) = \frac{p(G) x + p'(G)}{q(G) x + q'(G)}.
\end{equation}

\begin{mylem}\label{lem:relCylBad}
Fix a large parameter $\xi$. Let $\zeta = |\xi|^{\alpha}$ be a typical scale and, with $j(\zeta)$ as in \eqref{eq:typjScale}, suppose that $j_k \leq j(\zeta) < j_{k+1}$. Let $G = (\mathbf{a}_0, \ldots, \mathbf{a}_{j_k - 1}, b_k, \mathbf{a}_{j_k}, \ldots, \mathbf{a}_{j(\zeta)})$ be an admissible sequence. Let $I$ be a subinterval of $\mathbb{R}$ whose Lebesgue measure satisfies 
\begin{equation}\label{relCylBadAssn}
 |I| = |\xi|^{-1 + 2 \alpha \pm O(\epsilon)}.
\end{equation}
Then, provided that $|\xi|$ is sufficiently large, we have the estimate
\[\lambda_G^{\sharp}(I) \leq |I|^{\beta_F},\]
where
\[\beta_F := \frac{2 / \tau - 2 \alpha }{1 - 2 \alpha} - O(\epsilon).\]
\end{mylem}
\begin{proof}
Let $g$ denote the continued fraction map for this proof, so that $\lambda^{\sharp}(E) = \lambda(g^{-1}(E))$ and $\lambda_G^{\sharp}(E) \lambda_G(g^{-1}(E))$ for any Borel subset $E \subset \mathbb{R}$. Since $g$ is essentially bijective, we have that $\lambda(S) = \lambda^{\sharp}(g(S))$ and $\lambda_G(S) = \lambda_G^{\sharp}(g(S))$ for any set $S \subset \mathbb{N}^{\infty}$ in the cylinder $\sigma$-algebra.

Let $I$ be an interval with $| I |= |\xi|^{-1 + 2 \alpha \pm O(\epsilon)}$. Then the definition of $\lambda_G^{\sharp}$ guarantees that 
\begin{IEEEeqnarray*}{rCl}
& \lambda_G^{\sharp}(I) \\
= & \lambda_G(g^{-1}(I)) \\
= & \lambda(G^*)^{-1} \lambda(G \cdot g^{-1}(I)) \\
= & \lambda(G^*)^{-1} \lambda^{\sharp}(g(G \cdot g^{-1}(I)))
\end{IEEEeqnarray*}
Hence, to obtain an upper bound on $\lambda_G^{\sharp}$, it is enough to obtain a \textit{lower} bound on $\lambda(G^*)$ and an \textit{upper} bound on $\lambda^{\sharp}(g(G \cdot g^{-1}(I)))$.

The lower bound on $\lambda(G^*)^{-1}$ is given by Proposition \ref{cylinderlowerbound}. This lower bound states 
\[\lambda(G^*) \gtrsim |\cyl(G)|^{1 - O(\epsilon)}.\]
But we have from the basic properties of continued fractions that $|\cyl(G)| = K(G)^{-2 \pm O(\epsilon)}$, and we have from the definition of $j(\zeta)$ that $K(G) \leq |\xi|^{\alpha + O( \epsilon)}$. Hence
\begin{IEEEeqnarray*}{rCl}
\lambda(G^*) & \gtrsim & |\xi|^{-2\alpha - O(\epsilon)} .
\end{IEEEeqnarray*}

The upper bound on $\lambda^{\sharp}(g (G \cdot g^{-1}(I)))$ is given by the observation \eqref{contfracshift} and Proposition \ref{prop:unifBallCond}. Writing $\frac{p'}{q'}$ and $\frac{p}{q}$ for the final two convergents of the rational number $g(G)$, we have
\begin{IEEEeqnarray*}{Cl}
& \lambda^{\sharp}(g(G \cdot g^{-1}(I))) \\
\leq & |g(G \cdot g^{-1}(I))|^{\frac{2}{\tau} - O(\epsilon)} \\
\leq & \left|\frac{p I + p'}{q I + q'} \right|^{\frac{2}{\tau} - O(\epsilon)}.\\
\end{IEEEeqnarray*}
An estimate for $\left|\frac{p I + p'}{q I + q'} \right|$ is given by using the fact that for $1 \leq x \leq N + 1$, we have
\begin{IEEEeqnarray*}{Cl}
& \left| \frac{d}{dx} \frac{p x + p'}{qx + q'} \right| \\
= & \left|  \frac{q'p - p'q}{(qx + q')^2} \right| \\
\leq & \frac{1}{(qx + q')^2} \\
\leq & q^{-2}.
\end{IEEEeqnarray*}
As such, we see the interval $I$ is stretched by a factor of $\leq q^{-2}$ under the map $x \mapsto (px+p')/(qx+q')$ and we have from \eqref{relCylBadAssn} that
\begin{IEEEeqnarray*}{rCl}
\lambda^{\sharp}(g(G \cdot g^{-1}(I))) & \leq & (q^{-2} |I|)^{\frac{2}{\tau} - O(\epsilon)}. \\ 
& \leq & (|\xi|^{-2 \alpha + O(\epsilon)} |\xi|^{- 1 + 2 \alpha + O(\epsilon)})^{\frac{2}{\tau} + O(\epsilon)} \\
& \leq & |\xi|^{- \frac{2}{\tau} + O(\epsilon)}
\end{IEEEeqnarray*}
Combining the estimates on $\lambda(G^*)$ and $\lambda^{\sharp}(g(G \cdot g^{-1}(I)))$, we get
\[\lambda_G^{\sharp}(I) \lesssim |\xi|^{2 \alpha - \frac{2}{\tau} + O(\epsilon)}.\]
Using the assumption that $|I| \geq |\xi|^{-1 + 2 \alpha -  O(\epsilon)}$, we get
\[\lambda_G^{\sharp}(I) \lesssim |I|^{\frac{\frac{2}{\tau} - 2\alpha - O( \epsilon)}{1 - 2 \alpha -  O(\epsilon)}}\leq|I|^{\frac{\frac{2}{\tau} - 2 \alpha }{1 - 2 \alpha } -  O( \epsilon)}.\]
The desired bound is obtained by folding the implicit constant into an $|I|^{-\epsilon}$ term, which is possible if $\xi$ is large enough depending on $\epsilon$.
\end{proof}
\begin{mylem}\label{lem:relCylGood}
Fix a suitably large parameter $|\xi|$. Suppose that $|\xi|^{\alpha}$ is an exceptional scale, as defined in Section \ref{sec:endOfArg}, such that $(1-2\epsilon)j_k <\log|\xi|^{\alpha} < (\tau - 1+2\epsilon)j_{k}$. Set $\alpha'=(\tau - 1+10\epsilon)\alpha$ and $\zeta = |\xi|^{\alpha'}$. Then $\zeta$ is a typical scale. Furthermore, let $I \subset \mathbb{R}$ be an interval whose Lebesgue measure satisfies
\[ | I | = |\xi|^{-1 + 2 \alpha' \pm O(\epsilon)}.\]
 Then, for an admissible sequence $G=(\mathbf{a}_0, \ldots, \mathbf{a}_{j_k - 1}, b_k, \mathbf{a}_{j_k}, \ldots, \mathbf{a}_{j(\zeta)})$, we have that
\[\lambda_{G}(I)  \leq C |I|^{\beta_F},\]
where
\[\beta_F := 1-\epsilon.\]
\end{mylem}
\begin{myrmk}
The proof is wholly analogous to the proof of Lemma \ref{cylinderlemma} and is a simple consequence of our measure construction in Section \ref{Measure}, so is omitted. Roughly speaking, it can by understood by the fact that for cylinders $\cyl(G)$ of scale roughly $|\xi|^{-2\alpha'}$, it is the measures $\bar{\nu}$ which make substantial contributions to the size of $ \lambda (\cyl(G))$, because the next exceptional partial quotient is too far away.
\end{myrmk}

\section{Approximation of the relative measures}\label{sec:measAppr}
\setcounter{subsection}{1}
Suppose that $\zeta=|\xi|^\alpha$ is a typical scale and let $j(\zeta)$ be as in \eqref{eq:typjScale}. Recall the definition of $S(\zeta)$ appearing in \eqref{measuredecomp1}: the family of admissible sequences $G$ containing $j(\zeta)$ $\mathbf{a}_i$'s. We partitioned $S(\zeta)$ into a small number of equivalence classes $\mathcal{A}_{M_1,M_2}$ depending on the continuants $K(G)$ and $K'(G)$ in Definition \ref{def:AM1M2class}. 

We will establish the following lemma.
\begin{mylem}\label{classlemma}
Let $G \in \mathcal{A}_{M_1, M_2}$. Then we have the estimate
\[|\widehat{\lambda^\sharp}(\xi)| \leq O(|\xi|^{-\epsilon}) + \left|\sum_{M_1, M_2} \sum_{G \in \mathcal{A}_{M_1, M_2}} \lambda^\sharp(\cyl(G)) \int e \left( \xi  \frac{px + p'}{qx + q'} \right) \, d \lambda_{M_1, M_2}^\sharp(x) \right|.\]
\end{mylem}
\begin{proof}
Let $G \in \mathcal{A}_{M_1, M_2}$. Then elements of $\supp \, \lambda_{G}$ have the form
\[(\mathbf{a}_{j(\zeta) + 1}', \ldots, \mathbf{a}_{j_{k+1} - 1}', b_{k+1}', \ldots),\]
where each $\mathbf{a}_j$ is a $p$-vector of integers from $1$ to $N$, and $b_{k+1}$ is an integer satisfying the condition
\[\left(1 + \frac{\eta}{1000} \right)^{\gamma_{\eta}(G \cdot H)} \leq b_{k+1} < \left(1 + \frac{\eta}{1000} \right)^{\gamma_{\eta}(G\cdot H)},\]
where $H$ is the tuple
\[H := (\mathbf{a}_{j(\zeta) + 1}, \ldots, \mathbf{a}_{j_{k+1} - 1}).\]

The above tuple $H$ lies in the set $\{1, 2, \ldots, N\}^{p (j_{k + 1} - j(\zeta) - 1)}$. We will use $\mathcal{T}$ to denote the collection of these tuples. For $G\in A_{M_1,M_2}$, observe that the definition of $\lambda$ implies that for any $H \in \mathcal{T}$, $\lambda_{M_1, M_2}(H) = \lambda_G(H)$. Moreover, there is a fixed $\theta$ depending on $\zeta$ but not on $G$ or on $H$ such that for all $H\in\mathcal{T}$ with $H^* \cap \supp \lambda_G \neq \emptyset$, we have the bounds
\begin{equation}\label{thetadef}
\theta^{1 - \epsilon} < K(H) < \theta^{1 + \epsilon}.
\end{equation}
Recall the definition of the exponents $\gamma_\eta$ from the measure construction, \eqref{eq:gammaEtaDef}. Given an element $G \in \mathcal{A}_{M_1, M_2}$, define $\mathcal{T}_1(G)$ to be the set of sequences $H\in\mathcal{T}$ such that $\gamma_\eta(G \cdot H) = \gamma_{\eta}(G_{M,N} \cdot H)$, and define $\mathcal{T}_2(G)$ to be the set of sequences such that $\gamma_{\eta}(G \cdot H) \neq \gamma_{\eta}(G_{M,N} \cdot H)$. We make the following three claims:
\begin{myclm}\label{largetheta}
If $G \in \mathcal{A}_{M_1,M_2}$ and $\theta > |\xi|^{5}$, then 
\[\left|\int e \left(\xi  \frac{px + p'}{qx + q'} \right) \, d(\lambda_G^\sharp - \lambda_{M_1, M_2}^\sharp)(x) \right| \lesssim |\xi|^{-1}.\]
\end{myclm}
\begin{myclm}\label{samegamma}
If $G \in \mathcal{A}_{M_1, M_2}$, then 
\[\left|\int_{\cyl(\mathcal{T}_1)} e \left(\xi  \frac{px + p'}{qx + q'} \right) \, d(\lambda_G^\sharp - \lambda_{M_1, M_2}^\sharp)(x) \right| \lesssim |\xi|^{-1}.\]
\end{myclm}
\begin{myclm}\label{differentgamma}
If $\theta \leq |\xi|^{5}$, then, for $G\in \mathcal{A}_{M_1,M_2}$, the measures $\lambda_{G}^\sharp(\cyl(\mathcal{T}_2))$ and $\lambda_{M_1, M_2}^\sharp(\cyl(\mathcal{T}_2))$ are bounded above by $|\xi|^{-\epsilon}$.
\end{myclm}
We will now finish the proof of Lemma \ref{classlemma} assuming these three claims. The sum \eqref{measuredecomp1} can be decomposed by splitting the sum over $G$ into sums over the classes $\mathcal{A}_{M_1, M_2}$. This gives the decomposition
\begin{equation}\label{measuredecomp2}
\widehat{\lambda^\sharp}(\xi) = \sum_{M_1, M_2} \sum_{G \in \mathcal{A}_{M_1, M_2}} \lambda^\sharp(\cyl(G)) \int e \left(\xi  \frac{px + p'}{qx + q'} \right) \, d \lambda_G^\sharp(x).
\end{equation}

In the case in which $\theta > |\xi|^5$, Claim \ref{largetheta} shows that the measure $\lambda_G^\sharp$ can be replaced in each summand by $\lambda_{M_1, M_2}^\sharp$ with a $O(|\xi|^{-1})$ error term. Since the terms $\lambda^\sharp(\cyl(G))$ sum to $1$, replacing each $\lambda_G^\sharp$ by $\lambda_{M_1, M_2}^\sharp$ incurs only a $|\xi|^{-1}$ error term, proving the lemma in this case.

It only remains to handle the case in which $\theta \leq |\xi|^5$. We will approximate $d \lambda_G^\sharp$ by the measure $d \lambda_{M_1, M_2}^\sharp$. We will split the inner integral into an integral over $\cyl(\mathcal{T}_1)$ and an integral over $\cyl(\mathcal{T}_2)$. So we have $\widehat{\lambda^\sharp}(\xi) = \mathrm{I} + \mathrm{II} + \mathrm{III}$, where
\begin{IEEEeqnarray*}{rCl}
\mathrm{I} & = & \sum_{M_1, M_2} \sum_{G \in \mathcal{A}_{M_1, M_2}} \lambda^\sharp(\cyl(G)) \int_{\cyl(\mathcal{T}_1)} e \left(\xi \frac{px + p'}{qx + q'} \right) \, d (\lambda_G^\sharp - \lambda_{M_1, M_2}^\sharp)(x). \\
\mathrm{II} & = & \sum_{M_1, M_2} \sum_{G \in \mathcal{A}_{M_1, M_2}} \lambda^\sharp(\cyl(G)) \int_{\cyl(\mathcal{T}_2)} e \left(\xi \frac{px + p'}{qx + q'} \right) \, d (\lambda_G^\sharp - \lambda_{M_1, M_2}^\sharp)(x). \\
\mathrm{III} & = & \sum_{M_1, M_2} \sum_{G \in \mathcal{A}_{M_1, M_2}} \lambda^\sharp(\cyl(G)) \int e \left( \xi \frac{px + p'}{qx + q'} \right) \, d \lambda_{M_1, M_2}^\sharp(x).
\end{IEEEeqnarray*}
By Claim \ref{samegamma}, each integral appearing in $\mathrm{I}$ is bounded in magnitude by $O(|\xi|^{-1})$. Together with the fact that $\lambda^\sharp$ is a probability measure and thus $\sum_G \lambda^\sharp(\cyl(G)) = 1$, this proves that $\mathrm{I}$ is bounded above by $O(|\xi|^{-1})$.

By Claim \ref{differentgamma}, each integral occurring in $\mathrm{II}$ is bounded above by $O(|\xi|^{-\epsilon})$. Again using the fact that the terms $\lambda^\sharp(\cyl(G))$ sum to $1$, we see that $\mathrm{II}$ is bounded above by $O(|\xi|^{-\epsilon})$, proving the lemma.
\end{proof}
It remains to prove Claim \ref{largetheta}, Claim \ref{samegamma} and Claim \ref{differentgamma}.
\begin{proof}[Proof of Claim \ref{largetheta}]
Let $G \in \mathcal{A}_{M_1, M_2}$ and let $H \in \mathcal{T}$. We wish to estimate the integral
\begin{equation}\label{largethetaeqn}
\left| \int_{\cyl (H)} e \left(\xi \frac{p x + p'}{qx + q'} \right) d( \lambda_{G_{M_1, M_2}}^\sharp - \lambda_G^\sharp)(x) \right|.
\end{equation}
We have already seen that $\lambda_G^\sharp(\cyl(H)) = \lambda_{M_1, M_2}^\sharp(\cyl(H))$. Therefore, the integral \eqref{largethetaeqn} is bounded above by
\begin{IEEEeqnarray*}{Cl}
& 2 \lambda_G^\sharp(\cyl(H)) \sup_{x, y \in \cyl(H)} \left| e \left(\xi \frac{px + p'}{qx + q'}\right) - e \left(\xi \frac{py + p'}{qy + q'} \right) \right| \\
\lesssim & \lambda_G^\sharp(\cyl(H)) |\xi|\sup_{x, y \in \cyl(H)}  \left| \frac{px + p'}{qx + q'} - \frac{py + p'}{qy + q'} \right| \\
= & \lambda_G^\sharp(\cyl(H)) \sup_{x, y \in \cyl(H)} |\xi| \left| \frac{x - y}{(qx + q')(qy + q')} \right|.
\end{IEEEeqnarray*}
Moreover, the choice of $H$ guarantees that $\sup_{x, y \in \cyl(H)} |x - y| \lesssim \theta^{-2+2\epsilon}$. This gives the estimate
\begin{IEEEeqnarray*}{Cl}
& \left| \int_{\cyl(H)} e \left(\xi \frac{p x + p'}{qx + q'} \right) d( \lambda_{G_{M_1, M_2}}^\sharp -  \lambda_G^\sharp)(x) \right| \\
\lesssim & \lambda_G^\sharp(\cyl(H)) |\xi| \theta^{-2+2\epsilon} \\
\lesssim & \lambda_G^\sharp(\cyl (H)) |\xi|^{-8}.
\end{IEEEeqnarray*}
Summing over all $H \in \mathcal{T}$ and using the fact that the coefficients $\lambda_G^\sharp(\cyl(H))$ sum to one gives the desired estimate.
\end{proof}
\begin{proof}[Proof of Claim \ref{samegamma}]
Let $G \in \mathcal{A}_{M_1,M_2}$ and let $H \in \mathcal{T}_1(G)$. Then $\gamma_{\eta}(G \cdot H) = \gamma_{\eta}(G_{M_1,M_2} \cdot H)$. In particular, the definition of the measures $\lambda_G$ and $\lambda_{G_{M_1, M_2}}$ imply that for finite sequences $F$ of the form
\[F := (\mathbf{a}_{j(\zeta) + 1}, \ldots, \mathbf{a}_{j_{k+1} - 1}, b_{k+1}, \mathbf{a}_{j_{k+1}}, \ldots, \mathbf{a}_{j_{k+2} - 1}),\]
we have the equality $\lambda_G^\sharp(\cyl(F)) = \lambda_{G_{M_1, M_2}}^\sharp(\cyl(F))$. In particular, $\cyl(F)$ intersects the support of $\lambda_G$ if and only if it intersects the support of $\lambda_{G_{M_1, M_2}}$. 

We consider those $F$ as above for which $\cyl(F)$ intersects the support of $\lambda_G$. We work to estimate the integral
\begin{equation}\label{eq:cylFErrInt}\left|\int_{\cyl(F)} e \left(\xi \frac{px + p'}{qx + q'} \right) d( \lambda_{G_{M_1, M_2}}^\sharp -  \lambda_G^\sharp)(x) \right|,\end{equation}
where $\frac{p'}{q'}$ and $\frac{p}{q}$ are the final two convergents of the finite continued fraction whose partial quotients are given by $G$. 

This integral is easily seen to be bounded above by
\[2 \lambda_G^\sharp(\cyl(F)) \sup_{x,y \in \cyl(F)} \left| e \left(\xi \frac{px + p'}{qx + q'} \right) - e \left(\xi \frac{py + p'}{qy + q'} \right) \right|.\]

Suppose $x$ and $y$ lie in $\cyl(F)$. Then
\begin{IEEEeqnarray*}{rCl}
\left| e \left(\xi \frac{px + p'}{qx + q'} \right) - e \left(\xi \frac{py + p'}{qy + q'} \right) \right| & \lesssim & |\xi| \left| \frac{px + p'}{qx + q'} - \frac{py + p'}{qy + q'} \right| \\
& = & |\xi| \left| \frac{x - y}{(qx + q')(qy + q')} \right|
\end{IEEEeqnarray*}
If $j_{k+2}$ is chosen to be sufficiently large relative to $j_{k+1}$, then $|x - y| < |\xi|^{-100}$ whenever $x,y \in \cyl(F)$.  Hence we get
\[\left| e \left(\xi \frac{px + p'}{qx + q'} \right) - e \left(\xi \frac{py + p'}{qy + q'} \right) \right| \lesssim |\xi| \cdot q^{-2} \cdot |\xi|^{-100} \lesssim |\xi|^{-1}.\] 
Summing the resulting estimate for \eqref{eq:cylFErrInt} over all $F$ of the form $H \cdot F'$ where $H \in \mathcal{T}_1$ and using the fact that $G$ is a probability measure gives the desired estimate.
\end{proof}
\begin{proof}[Proof of Claim \ref{differentgamma}]
Let $G \in \mathcal{A}_{M_1, M_2}$ and suppose $H \in \mathcal{T}_2(G)$. Let $\frac{p'}{q'}$ and $\frac{p}{q}$ be the last two convergents of the finite continued fraction with partial quotients given by $G$, and let $\frac{\tilde p'}{\tilde q'}$ and $\frac{\tilde p}{\tilde q}$ be the final two partial quotients in the continued fraction expansion of $G_{M_1, M_2}$. Let $\frac{p^*}{q^*}$ be the rational number whose finite continued fraction expansion is given by $H$. Since $H \in \mathcal{T}_2(G)$, we have for some integer $\gamma \lesssim \log |\xi|$ that one of the following two inequalities holds:
\begin{equation}\label{badineq1}
q' p^* +  q q^* < \left(1 + \frac{\eta}{1000} \right)^{\gamma} \leq \tilde q' p^* + \tilde q q^*
\end{equation}
or
\begin{equation}\label{badineq2}
\tilde q' p^* +  \tilde q q^* < \left(1 + \frac{\eta}{1000} \right)^{\gamma} \leq q' p^* + q q^*.
\end{equation}
Without loss of generality, suppose \eqref{badineq1} holds for a given $\gamma$. Since $\tilde q' = q' + O(|\xi|^{\alpha - 200 \epsilon})$ and $\tilde q = q + O(|\xi|^{\alpha - 200 \epsilon})$, we must then have that
\[|(q' - \tilde q') p^*| + |(q - \tilde q)q^*| \leq |\xi|^{\alpha - 200 \epsilon}  \max(p^*, q^*).\]
Hence \eqref{badineq1} can only occur if
\[\left|q' p^* + q q^* - \left(1 + \frac{\eta}{1000} \right)^{\gamma} \right| \leq |\xi|^{\alpha - 200 \epsilon} \max(p^*, q^*) .\]
Recall that $q^*$ satisfies $\theta^{1 - \epsilon} < q^* < \theta^{1 + \epsilon}$ for an appropriate quantity $\theta$ not depending on $H$. Since $p^* \leq (N + 1) q^*$, we have the estimate
\begin{IEEEeqnarray*}{rCl}
\left| p^* + \frac{q}{q'} q^* - \frac{1}{q'} \left(1 + \frac{\eta}{1000} \right)^{\gamma} \right| & \leq & |q'|^{-1} |\xi|^{\alpha - 200 \epsilon} \max(p^*, q^*) \\
& \leq & |\xi|^{- 199 \epsilon} \max(p^*, q^*).
\end{IEEEeqnarray*}
Hence for a fixed $q^*$ and $\gamma$, there are only $O(|\xi|^{-199 \epsilon}) \cdot q^* $ choices available for $p^*$. Moreover, it is easy to verify that $q^* \leq p^* \leq (N+1) q^* \leq \theta^{1 + 2 \epsilon}$, provided that $\theta$ is sufficiently large. Hence the total number of choices for the pair $(p^*, q^*)$---that is to say, the total number of choices for $H$---is bounded above by $O_N(\theta^{2 + 3\epsilon} |\xi|^{-199 \epsilon}) \leq O(\theta^{2} |\xi|^{- 184 \epsilon}).$

We now take advantage of the relative ball condition satisfied by the measures $\lambda_G^\sharp$ and $\lambda_{M,N}^\sharp$. By an argument analogous to the proof of Lemma \ref{cylinderlemma}, we have that for each $H$, $\lambda_G^\sharp$ satisfies an inequality of the form
\[\lambda_G^\sharp(\cyl(H)) \lesssim |\cyl(H)|^{1 - \epsilon} \lesssim |\theta|^{-2 + 4 \epsilon} \lesssim |\theta|^{-2} |\xi|^{20\epsilon}.\]
Since there are at most $O(\theta^2 |\xi|^{-184 \epsilon})$ choices for $H$ for each $\gamma$, and at most $O(\log |\xi|) \leq O(|\xi|^{\epsilon})$ choices for $\gamma$, we therefore have that 
\[\lambda_G^\sharp(\cyl(\mathcal{T}_2)) \lesssim |\xi|^{-160 \epsilon},\]
as desired.
\end{proof}

\appendix

\section{Appendix: Estimates of Kaufman}\label{appendixKaufmanEsts}
\setcounter{subsection}{1}
We will summarize the proofs of the key estimates of Kaufman, characterised by Lemmas \ref{QRlem}, \ref{Mlemma}, and \ref{m2lemma}. These proofs essentially appear in Kaufman \cite{Kaufman80}, but are included here for completeness.

Recall the decomposition from Section \ref{sec:combEst}. In our case, we consider
\[F(x)=f_{\xi}(x)=\sum_{G\in \mathcal{A}_{M_1,M_2}} \lambda^\sharp(\cyl(G)) e\left(-\xi\dfrac{p(G)x+p'(G)}{q(G)x+q'(G)}\right)\]
and $\|f_{\xi}\|_{L^2([1,N+1])}^2$ can be treated as a classical oscillatory integral. 

\begin{proof}[Proof of Lemma \ref{Mlemma}]
Observe that, by construction, $|\xi|^{\alpha - \epsilon} < q(G) < |\xi|^{\alpha + \epsilon}$. By the triangle inequality, we obtain 
\[\left|f_{\xi}'(x)\right|\leq\sum_{G\in \mathcal{A}_{M_1,M_2}} \lambda^\sharp(\cyl(G)) \left|\frac{d}{dx}e\left(-\xi\dfrac{p(G)x+p'(G)}{q(G)x+q'(G)}\right)\right|.\]
For each of the cofactors of $\lambda^\sharp(\cyl(G))$, we find, using $pq' - p'q = \pm 1$, for $1\leq x \leq N+1$,
\[\left|2\pi\xi\dfrac{p(G)(q(G)x+q'(G))-(p(G)x+p'(G))q(G)}{(q(G)x+q'(G))^2}\right|\]
\[=\left|2\pi\xi\dfrac{1}{(q(G)x+q'(G))^2}\right|\]
\[\leq\left|\xi\right|^{1-2\alpha+3\epsilon},\]
so taking the sum and using the fact $\lambda^\sharp$ is a probability measure gives the desired bound.
\end{proof}
In order to prove Lemma \ref{m2lemma}, Kaufman establishes two van der Corput-type inequalities. The first is useful for nonstationary phases $f$.
\begin{mylem}\label{nonstationaryvdc}
Let $f$ be a real-valued $C^2$ function on the interval $[a,b]$, and suppose $A$ and $B$ are positive numbers such that $f' \geq A$ or $f' \leq -A$,  and $|f''| \leq B$ on $[a,b]$. Then we have the estimate
\[\left| \int e(f(x)) \, dx \right| \leq \frac{1}{\pi} A^{-1} + \frac{b - a}{2 \pi} A^{-2} B.\]
\end{mylem}
\begin{proof}
Integrate by parts, using $u(x) = \frac{1}{f'(x)}$ and $v'(x) = f'(x) e(f(x))$.  Then the integration by parts formula gives
\begin{IEEEeqnarray*}{rCl}
\left| \int_a^b e(f(x)) \, dx \right| & \leq & \left| \frac{1}{2 \pi} \frac{e(f(x))}{f'(x)} \right|_a^b + \frac{1}{2 \pi} \left| \int_a^b e(f(x)) \frac{f''(x)}{(f'(x))^2} dx\right| \\
& \leq & \frac{1}{\pi} A^{-1}+ \frac{b - a}{2 \pi} A^{-2} B,
\end{IEEEeqnarray*}
as desired.
\end{proof}
The second van der Corput type lemma concerns the situation in which the derivative of the phase is the product of a linear function and a nonvanishing function. As the phase is stationary, the estimate will be weaker than that of Lemma \ref{nonstationaryvdc}.
\begin{mylem}\label{stationaryvdc}
Suppose that $f'(x) = (C_1x + C_2) g(x)$ on some compact interval $[a,b]$, where $g$ is a $C^1$ function satisfying the estimates $|g(x)| \geq A$, $|g'(x)| \leq B$ for $x \in [a,b]$, where $B > A$. Then there exists an absolute constant $K$ such that we have the oscillatory integral estimate
\[\left| \int_a^b e(f(x)) \, dx \right|< K (1 + (b-a)) B A^{-3/2} |C_1|^{-1/2}.\]
\end{mylem}
\begin{proof}
Write 
\[[a,b] = I_L \cup J \cup I_R\]
where $J$ is the interval $[-\frac{C_2}{C_1} - |C_1 A|^{-1/2}, -\frac{C_2}{C_1} + |C_1 A|^{-1/2}]$.
The complement of $J$ in $[a,b]$ can be written as $I_L \cup I_R$, where each of $I_L$ and $I_R$ is an interval (possibly empty) containing the points in $[a,b]$ to the left and right of $J$, respectively. The length of $J$ is at most $2 |A C_1|^{-1/2}$. Since $B > A$, we have that this is bounded above by $2A^{-3/2} B|C_1|^{-1/2}$.

We will estimate the integral over $I_L$; the integral over $I_R$ is estimated in a similar manner. Before estimating these integrals, we need some estimates on the derivatives of $f$.

On each of $I_L$ and $I_R$ we have $|C_1x + C_2| \geq C_1^{1/2} A^{-1/2}$. Therefore, we have the estimate
\[|f'(x)| = |C_1 x + C_2| |g(x)| \geq |C_1|^{1/2} A^{-1/2} A = |A C_1|^{1/2}.\]
Before estimating the integral, we observe from the product rule that
\[f''(x) = C_1 g(x) + (C_1 x + C_2) g'(x).\]
Therefore,
\begin{IEEEeqnarray*}{rCl}
\left|\frac{f''(x)}{f(x)^2} \right| & = & \left| \frac{C_1 g(x) + (C_1 x + C_2)g'(x)}{(C_1x + C_2)^2 g(x)^2} \right| \\
& \leq & |C_1| \left| \frac{1}{g(x) (C_1 x + C_2)^2} \right| + \left|\frac{(C_1 x + C_2) g'(x)}{(C_1 x + C_2)^2 g(x)^2} \right| \\
& \leq & \frac{|C_1| A^{-1}}{(C_1 x + C_2)^2} + \frac{B A^{-2}}{|C_1 x + C_2|}.
\end{IEEEeqnarray*}
Now, we are ready to integrate by parts as in Lemma \ref{nonstationaryvdc}. On $I_L$, this integral is
\begin{IEEEeqnarray*}{rCl}
\left| \int_{a_L}^{b_L} e (f(x)) \, dx \right|  & \leq & \left| \frac{1}{2 \pi} \frac{e(f(x))}{f'(x)} \right|_{a_L}^{b_L} + \frac{1}{2 \pi} \left| \int_{a_L}^{b_L} e(f(x)) \frac{f''(x)}{(f'(x))^2} dx\right|\\
& \leq & \frac{b_L - a_L}{2 \pi |A C_1|^{1/2}} + \frac{1}{2 \pi} \int_{a_L}^{b_L}\left( \frac{|C_1| A^{-1}}{(C_1 x + C_2)^2} + \frac{A^{-2}B}{|C_1x + C_2|} \right)dx 
\end{IEEEeqnarray*}
It remains to estimate both integrals
\[\mathrm{I} := \frac{1}{2 \pi} \int_{a_L}^{b_L} \frac{|C_1| A^{-1}}{(C_1 x + C_2)^2} \, dx\]
and
\[\mathrm{II} := \frac{1}{2 \pi} \int_{a_L}^{b_L} \frac{A^{-2} B}{C_1x + C_2}dx.\]
First, we will estimate the integral $\mathrm{I}$. Keeping in mind that $|C_1 x + C_2| \geq |C_1|^{1/2} A^{-1/2}$ on the interval $[a_L, b_L]$, we have
\begin{IEEEeqnarray*}{Cl}
& \frac{1}{2 \pi} \int_{a_L}^{b_L} \frac{|C_1| A^{-1}}{(C_1 x + C_2)^2} \, dx  \\
= & \frac{1}{2 \pi} \int_{x = a_L}^{x = b_L} \frac{|C_1| A^{-1}}{u^2} \cdot  \frac{1}{|C_1|} \, du \\
= & \frac{1}{2 \pi} A^{-1} \int_{x = a_L}^{x = b_L} u^{-2} \, du \\
= & \left| \frac{1}{2 \pi} A^{-1} u^{-1} \right|_{x = a_L}^{x = b_L} \\
= & \left| \frac{1}{2 \pi} A^{-1} (C_1 x + C_2)^{-1} \right|_{a_L}^{b_L} \\
\leq & 2 \cdot \frac{1}{2 \pi} A^{-1} |C_1|^{-1/2} A^{1/2}
\lesssim A^{-3/2} B |C_1|^{-1/2}.
\end{IEEEeqnarray*}
As for the remaining term, we again use the fact that $|C_1x + C_2| \geq A^{-1/2} |C_1|^{1/2}$ on $[a_L, b_L]$.
\begin{IEEEeqnarray*}{rCl}
\mathrm{II} & = &  \frac{1}{2 \pi} \int_{a_L}^{b_L} \frac{A^{-2} B}{|C_1 x + C_2|} dx\\
& \leq & \frac{1}{2 \pi} \frac{A^{-2} B}{A^{-1/2} |C_1|^{1/2}} (b_L - a_L) \\
& \leq & \frac{1}{2 \pi} A^{-3/2} B |C_1|^{-1/2}(b_L - a_L).
\end{IEEEeqnarray*}
Applying the same estimate to the integral over $I_R$ and adding all of the terms together gives the desired bound.
\end{proof}

\begin{proof}[Proof of Lemma \ref{QRlem}]
Let $r > 0$. Decompose the interval $[a,b]$ into a collection $\mathcal{I}$ of smaller intervals $I$ with length $\frac{r}{M}$, with one possible shorter interval $I^*$ that is not counted as being in $\mathcal{I}$. Observe that if $|F(x_0)| \geq 2r$ at some point in any such interval $I$, then in fact the bound $|F'(x)| \leq M$ implies that $|F(x)| \geq r$ on all of $I$. Let $\mathcal{I}_{\text{Large}}$ be the collection of intervals $I$ containing some point $x_0$ with $F(x_0) \geq 2r$ and $\mathcal{I}_{\text{Small}}$ denote $\mathcal{I} \setminus \mathcal{I}_{\text{Large}}$. Define
\begin{IEEEeqnarray*}{rCl}
E_{\text{Large}} & := & I^* \cup \bigcup_{I \in \mathcal{I}_{\text{Large}}} I \\
E_{\text{Small}} & := & \bigcup_{I \in \mathcal{I}_{\text{Small}}} I 
\end{IEEEeqnarray*}
We count $\mathcal{I}_{\text{Large}}$ using a pigeonhole principle argument: 
\[m_2 \geq \#\mathcal{I}_{\text{Large}} \cdot r^2 \cdot \frac{r}{M}\]
and thus
\begin{equation}
\#\mathcal{I}_{\text{Large}} \leq m_2 M r^{-3}.
\end{equation}
Since each $I \in \mathcal{I}_{\text{Large}}$ satisfies $\mu^\sharp(I) \leq |I|^{\beta} = \left(\frac{r}{M} \right)^{\beta}$, we have that $\mu^\sharp(E_{\text{Large}}) \leq \left(\frac{r}{M} \right)^{\beta} (1 + m_2 M r^{-3})$.
On the set $E_{\text{Large}}$, we apply the estimate $|F(x)| \leq 1$; this gives the estimate
\[\int_{E_{\text{Large}}} |F(x)| \, d \mu^\sharp(x) \leq \left(\frac{r}{M} \right)^{\beta} (1 + m_2 M r^{-3}).\]
On the set $E_{\text{Small}}$, we apply the estimate $|F(x)| \leq 2r$. Since $\mu^\sharp$ is a probability measure, we have
\[\int_{E_{\text{Small}}}|F(x)| d \mu^\sharp(x) \leq 2r.\]
So in total, we have
\[\int_a^b |F(x)| d \mu^\sharp(x) \leq 2r + \left(\frac{r}{M} \right)^{\beta} (1 + m_2 M r^{-3}),\]
as desired.
\end{proof}

Now, we turn our attention to Lemma \ref{m2lemma}. The analysis closely follows Kaufman.
\begin{proof}It is seen directly that
\[\|f_{\xi}\|_{L^2([1,N+1])}^2=\sum_{G\in \mathcal{A}_{M_1,M_2}} \int_{1}^{N+1}\sum_{\tilde{G}\in \mathcal{A}_{M_1,M_2}} \lambda^\sharp(\cyl(G))\lambda^\sharp(\cyl(\tilde{G})) e\left(\Phi_{(G,\tilde{G})}(x)\right)dx,\]
where \[\Phi_{(G,\tilde{G})}(x)=-\left(\dfrac{p(G)x+p'(G)}{q(G)x+q'(G)}-\dfrac{p(\tilde{G})x+p'(\tilde{G})}{q(\tilde{G})x+q'(\tilde{G})}\right).\]
Considering, $\Phi'_{(G,\tilde{G})}$, it is easy to verify that, up to a factor of $\pm1$,
\[\Phi_{(G,\tilde{G})}'(x)=\xi\left(\dfrac{1}{(q(G)x+q'(G))^2}-\dfrac{1}{(q(\tilde{G})x+q'(\tilde{G}))^2}\right)\]
\[=g(x)(C_1x+C_2)\]
where $C_1=\xi(q(\tilde{G})-q(G))$, $C_2=\xi(q'(\tilde{G})-q'(G))$ and
\[g(x)=\dfrac{1}{(q(G)x+q'(G))^2(q(\tilde{G})x+q'(\tilde{G}))^2}\left(((q(\tilde{G})+q(G))x+(q'(\tilde{G})+q'(G)))\right).\]
When $q(G)=q(\tilde{G})$ and $q'(G)\neq q'(\tilde{G})$, it is clear that the phase is nonstationary. We classify such pairs $(G,\tilde{G})\in\mathcal{C}_1$. When $q(G)\neq q(\tilde{G})$, it is clear we have at most one critical point in the region of integration, we classify such pairs $(G,\tilde{G})\in\mathcal{C}_2$. Finally, where $q(G)=q(\tilde{G})$ and $q'(G)= q'(\tilde{G})$, the phase is constant. We classify such pairs  $(G,\tilde{G})\in\mathcal{C}_3$.

For the pairs $(G,\tilde{G})\in\mathcal{C}_1$, the phase is nonstationary. In this case, $\Phi'_{(G,\tilde{G})}$ has the form
\[\Phi'(x) = \xi  \frac{((q(G)+ q(\tilde{G}))x + q'(G) + q'(\tilde{G}))(q'(G) - q'(\tilde{G}))}{(q(G) x +q'(G))^2 (q(\tilde{G}) x + q'(\tilde{G}))^2}\]
Keeping in mind that $|\xi|^{\alpha - \epsilon} < q(G), q'(G), q(\tilde{G}), q'(\tilde{G}) < |\xi|^{\alpha + \epsilon}$ and that $1 \leq x \leq N + 1$, we have that $\left|\Phi'(x)\right|$ is bounded below by $C_N |q'(G) - q'(\tilde{G})| |\xi| \cdot |\xi|^{\alpha - \epsilon} \cdot |\xi|^{-4 \alpha - 4 \epsilon} \geq |q'(G) - q'(\tilde{G})| |\xi|^{1 - 3 \alpha +O( \epsilon)}$, provided that $|\xi|$ is sufficiently large depending on $N$. We compute the second-derivative of $\Phi$. This second-derivative will be of the form
\[\Phi''(x) = \xi \frac{F_1(q(G), q'(G), q(\tilde{G}), q'(\tilde{G}), x) (q'(G) - q'(\tilde{G}))}{F_2(q(G), q(\tilde{G}), x)}\]
where $F_1$ is a polynomial of degree five in $q(G), q'(G), q(\tilde{G}),$ and $q'(\tilde{G})$, and $F_2$ is a polynomial of degree $8$ in $q(G), q'(G)q(\tilde{G}), q(\tilde{G}),$ and $q'(\tilde{G})$. Hence, for large enough $|\xi|$, we have the estimate 
\[\left| \Phi''(x)\right| \leq C_N |q'(G) - q'(\tilde{G})| \cdot |\xi| \cdot |\xi|^{5 \alpha + 5 \epsilon} \cdot |\xi|^{-8 \alpha + 8 \epsilon} \leq |q'(G) - q'(\tilde{G})| |\xi|^{1 -3 \alpha + O(\epsilon)}.\]
Hence we can apply Lemma \ref{nonstationaryvdc} with the parameters
\begin{IEEEeqnarray*}{rCl}
A & = & |q'(G) - q'(\tilde{G})| |\xi|^{1 - 3 \alpha + O(\epsilon)} \\ 
B & = & |q'(G)- q'(\tilde{G})| |\xi|^{1 - 3 \alpha + O(\epsilon)}.
\end{IEEEeqnarray*}
This yields the estimate
\begin{IEEEeqnarray*}{Cl}
& \sum_{(G,\tilde{G})\in\mathcal{C}_1} \lambda^\sharp(\cyl(G))\lambda^\sharp(\cyl(\tilde{G}))\left|\int_{1}^{N+1} e\left(\Phi_{(G,\tilde{G})}(x)\right)dx\right| \\
\leq & \sum_{(G,\tilde{G})\in\mathcal{C}_1} \lambda^\sharp(\cyl(G))\lambda^\sharp(\cyl(\tilde{G}))|q'(G) - q'(\tilde{G})|^{-1} |\xi|^{-1 + 3 \alpha + O(\epsilon)} \leq |\xi|^{-1 + 3 \alpha + O(\epsilon)}.
\end{IEEEeqnarray*}

For the pairs $(G,\tilde{G})\in\mathcal{C}_2$, we may have a critical point for the phase. We have that $g\in C^1$ and $|g|\geq \frac{|\xi|^{\alpha-4\epsilon}}{100|\xi|^{4(\alpha+\epsilon)}}\geq|\xi|^{-3\alpha+ O(\epsilon)}=A$. We also have that 
\[|g'(x)|\leq C_N\left(\dfrac{|\xi|^{5(\alpha+\epsilon)}}{|\xi|^{8(\alpha-\epsilon)}}+\dfrac{|\xi|^{2(\alpha+\epsilon)}}{|\xi|^{5(\alpha-\epsilon)}}\right)\leq |\xi|^{-3\alpha+ O(\epsilon)}=B.\]
We can apply Kaufman's Lemma \ref{stationaryvdc} to obtain the bound
\begin{IEEEeqnarray*}{Cl}
& \sum_{(G,\tilde{G})\in\mathcal{C}_2} \lambda^\sharp(\cyl(G))\lambda^\sharp(\cyl(\tilde{G}))\left|\int_{1}^{N+1} e\left(\Phi_{(G,\tilde{G})}(t)\right)dt\right| \\
\leq & \sum_{(G,\tilde{G})\in\mathcal{C}_2} \lambda^\sharp(\cyl(G))\lambda^\sharp(\cyl(\tilde{G}))BA^{-\frac{3}{2}}|C_1|^{-\frac{1}{2}} \\
\leq & \sum_{(G,\tilde{G})\in\mathcal{C}_2} \lambda^\sharp(\cyl(G))\lambda^\sharp(\cyl(\tilde{G}))BA^{-\frac{3}{2}}|\xi|^{-\frac{1}{2}} \\
\leq & C_N''|\xi|^{\frac{3\alpha-1}{2} + O(\epsilon)}.
\end{IEEEeqnarray*}

For those pairs $(G,\tilde{G})\in\mathcal{C}_3$, the phase is constant. Here, we make use of the geometric information we have about the measure. By Lemma \ref{cylinderlemma}, we know that $\lambda^\sharp$ satisfies a ball condition on cylinders with exponent $\frac{\tau }{2 \tau - 2}+ O(\epsilon)$. If $y$ is a real number with $p(\tilde{G})/q(\tilde{G})$ as a convergent, then $\left|y-p(\tilde{G})/q(\tilde{G})\right|\leq q(\tilde{G})^{-2}\leq |\xi|^{-2\alpha(1-\epsilon)}$ so $\cyl(G)$ has diameter $\leq 2 |\xi|^{-2\alpha(1-\epsilon)}$. For $(G,\tilde{G})\in\mathcal{C}_3$, we have that $q(G)=q(\tilde{G})$ and $q'(G)= q'(\tilde{G})$. For a given $G$, we consider how many pairs $(G,\tilde{G})\in\mathcal{C}_3$. Observe that except for the integer term $\tilde{a}_0$, all of $\tilde{G}$ is entirely determined by $q(\tilde{G})$ and $q'(\tilde{G})$, since it is possible to use the Euclidean algorithm to determine the previous denominators and hence the partial quotients. Therefore, for a fixed $G$, there are at most $N$ possible choices for $\tilde{G}$. Thus we find
\begin{IEEEeqnarray*}{Cl}
& \left|\sum_{(G,\tilde{G})\in\mathcal{C}_3} \lambda^\sharp(\cyl(G))\lambda^\sharp(\cyl(\tilde{G}))\int_{1}^{N+1} e\left(\Phi_{(G,\tilde{G})}(x)\right)dx\right| \\
\leq & N\sum_{G} \sum_{\tilde{G}|\,(G,\tilde{G})\in\mathcal{C}_3}\lambda^\sharp(\cyl(G))\lambda^\sharp(\cyl(\tilde{G})) \\
\leq & N\sum_{G} \sum_{\tilde{G}|\,(G,\tilde{G})\in\mathcal{C}_3}\lambda^\sharp(\cyl(G))(2|\xi|^{-2\alpha(1-\epsilon)})^{\left(\frac{\tau }{2 \tau - 2}\right)+ O(\epsilon)} \\
\leq & N^2|\xi|^{-2\alpha\left(\frac{\tau}{2 \tau - 2}\right)+ O(\epsilon)}\sum_{G} \lambda^\sharp(\cyl(G)).\\
= & N^2|\xi|^{-2\alpha\left(\frac{\tau}{2 \tau - 2}\right)+ O(\epsilon)}.
\end{IEEEeqnarray*}

Summing the estimates, we find 
that
\[\|f_{\xi}\|_{L^2([1,N+1])}^2\leq C_N \left(|\xi|^{\frac{3\alpha-1}{2}+ O(\epsilon)}+|\xi|^{-\alpha\left(\frac{\tau }{ \tau - 1}\right) + O(\epsilon)}\right).\]
\end{proof}
\section*{Acknowledgements}
The authors would like to thank Sanju Velani and Evgeniy Zorin, without whom this project would not have been possible.
\bibliographystyle{myplain}
\bibliography{Diophantine_Exact_Order_Small_Tau}
\end{document}